\documentclass[11pt]{amsart}
\usepackage{latexsym}   % use all LaTeX fonts
\usepackage{amssymb}    % use all AMS fonts
\usepackage{amsmath}    % include some AMS-LaTeX functionality
\usepackage{amsthm}     % AMS style theorem and proof environments
\usepackage{mathrsfs}   % special calligraphic characters like \mathscr{S}
\usepackage{euscript}   % euscript package
\usepackage{epsfig}     % epsfig package, encap postscript
\usepackage[all]{xy}    % Xy-Pic for diagrams
\newcommand{\minus}{\setminus}
\newcommand{\cross}{\times}

\newcommand{\comp}{\circ}

\newcommand{\ra}{\rightarrow}
\newcommand{\R}{\mathbb{R}}

\renewcommand{\epsilon}{\varepsilon}
\newcommand{\nin}{\not\in}
\renewcommand{\qed}{\square}
\renewcommand{\phi}{\varphi}
\newcommand{\Ha}{\mathscr{H}}
\newcommand{\ip}{<\cdot,\cdot>}
\newcommand{\pf}{\noindent {\em Proof: }}
\newtheorem{Pro}{Proposition}[section]
\newtheorem{Lem}[Pro]{Lemma}

\newtheorem{Thm}[Pro]{Theorem}
\newtheorem{MThm}{Theorem}

\newtheorem{Def}[Pro]{Definition}
\newtheorem{Example}[Pro]{Example}

\newtheorem{Cor}[Pro]{Corollary}

\newcommand{\hook}{\lfloor}
\newcommand{\cc}{d_{cc}}

\newcommand{\f}{\frac}
\def\lra{\longrightarrow}
\newcommand{\n}{\nabla}
\newcommand{\ov}{\overline}
\def\sbq{\subseteq}

\begin{document}

\title{A notion of rectifiability modeled on Carnot groups}
\author{Scott D. Pauls}
%\date{\today}
\address{Dartmouth College, Hanover, NH 03755}
\email{scott.pauls@dartmouth.edu}
\thanks{The author was partially supported by NSF grants DMS-9971563
  and DMS-0306752}
\begin{abstract}
We introduce a notion of rectifiability modeled on Carnot groups.
Precisely, we say that a subset $E$ of a Carnot group $M$ and $N$ is a
subgroup of $M$, we say $E$ is $N$-rectifiable if it
is the Lipschitz image of a positive measure subset of $N$.  First, we
discuss the implications of $N$-rectifiability, where
$N$ is a Carnot group (not merely a subgroup of a Carnot group), which include
$N$-approximability and the existence of approximate
tangent cones isometric to $N$ almost everywhere in $E$.  Second, we
prove that, under a stronger condition concerning the existence of
approximate tangent cones isomorphic to $N$ almost everywhere in a set
$E$, that $E$ is 
$N$-rectifiable. Third, we
investigate the rectifiability properties of level sets of $C^1_N$
functions, $f:N \ra \R$, where $N$ is a Carnot group.  We show that
for almost every $t \in \R$ and almost every noncharacteristic $x \in f^{-1}(t)$, there
exists a subgroup $T_x$ of $H$  and $r >0$ so that $f^{-1}(t) \cap
B_H(x,r)$ is $T_x$-approximable at $x$ and an approximate tangent cone
isomorphic to $T_x$ at $x$.  
\end{abstract}

\maketitle
\section{Introduction}
The notion of rectifiability is central to the study of geometric measure theory, allowing for the proof of classical
geometric properties in a much more general setting.  In recent years,
there has been significant interest and progress in the study of
rectifiable sets not only in Euclidean space but in more general
metric spaces as well (see, for example,  \cite{Amb2}, \cite{AmbKir},
\cite{AmbKir2}, \cite{Cheeger}, \cite{brokendreams}, \cite{FSSC}, \cite{FSSC2},
\cite{GarNh}, \cite{Kir}, \cite{Magnani}, \cite{Pansu}, \cite{White}.
Also, see the extensive bibliographies in \cite{GarNh} and \cite{AmbKir}).
  In attempting to use the techniques of geometric
measure theory to investigate the properties of general metric spaces,
one quickly encounters a major difficulty:  there may not be any
rectifiable subsets or the set of rectifiable subsets may be too small
to reveal any significant geometry.  Here, we consider a subset of a
metric space to be rectifiable if it can be realized as the Lipschitz
image of a piece of Euclidean space.  Thus, to have any hope of
transporting the techniques of Euclidean geometric measure theory to
metric spaces, we need a more general notion of rectifiable sets which
are modeled on a wider class of metric spaces than simply Euclidean ones.  In
this paper, we investigate a special situation where much of the
standard rectifiable theory carries over, but reveals some of the
complications inherent in this endeavor.  

We will restrict ourselves to the investigation of the so-called
Carnot groups - connected, simply connected, graded nilpotent Lie
groups equipped with a left-invariant Carnot-Carath\'eodory metric
(see below for precise definitions).  Carnot groups arise in a variety
of situations:  in the asymptotic geometry of manifolds of negative
curvature, in optimal control theory, in the local geometry of
equiregular Carnot-Carath\'eodory manifolds, in CR geometry, in the study of
subelliptic partial differential equations, and in many other areas.
In addition, Carnot groups are 
analytically very similar to Euclidean spaces (which are themselves particular
examples of Carnot groups) in that they possess translations and
dilations that respect the metric.  Because of these similarities to
Euclidean spaces, Carnot groups form a good test class for the
generalization of the notion of rectifiability.

Further motivation for considering rectifiability modeled on Carnot
groups arises in the study of the local geometry of
Carnot-Carath\'eodory manifolds.  A Carnot-Carath\'eodory
manifold is a quadruple $(X,\mathcal{S},\ip,d_X)$ where $X$ is a
smooth manifold, $\mathcal{S}$ is a subbundle of the tangent bundle,
$\ip$ is a smoothly varying inner product on the fibers of
$\mathcal{S}$ and $d_X$ is the path metric formed by taking the
infimum of lengths of paths (calculated with respect to $\ip$) among
absolutely continuous paths tangent a.e. to $\mathcal{S}$.  Further, letting $\mathfrak{d}_i(x)$ be the
dimension of the
subspace of $T_xX$ spanned by all the commutators of order less than
or equal to $i$, we say that $(X,\mathcal{S},\ip,d_X)$ is equiregular if the
vector $(\mathfrak{d}_0(x),\mathfrak{d}_1(x),...)$ is locally constant.  Carnot groups are very special
examples of equiregular Carnot-Carath\'eodory spaces given by the
quadruple $(G, \mathcal{V},\ip,d_G)$ where $\mathcal{V}$ is the bottom
level of the 
grading, thought of as a left invariant subbundle of $TG$, and $\ip$
is a left invariant inner product on $\mathcal{V}$.  For general
equiregular Carnot-Carath\'eodory spaces, it is known
(\cite{Mitchell}, see also \cite{Bell:1996}) that the tangent
cone to $(X,\mathcal{S},\ip,d_X)$ is isometric to a Carnot group.
However, the tangent group may vary from
point to point (quite badly!) - see, for example, \cite{Bell:1996} and
\cite{Varchenko}.  Naturally, this makes an investigation of the local
geometry of Carnot-Carath\'eodory manifolds much more difficult than,
for example, that of Riemannian manifolds.  As we shall see, the
framework of rectifiability modeled on Carnot groups identifies a
subclass of CC manifolds which possess uniform local behavior.

Let $(M,d)$ be a Carnot group and let $N$ be
a subgroup of $M$ equipped with a Carnot metric $d'$.  We define the following
generalization of rectifiability:  a subset $S \subset M$ is
$N$-rectifiable if it is the Lipschitz (with respect to $d$ and $d'$)
image of a positive measure subset of $N$.  This clearly generalizes the
standard notion of rectifiability where $M$ and $N$ are replaced by
Euclidean spaces of the appropriate dimension.  Using this
terminology, we will refer to the standard theory of rectifiability as
$\R^k$-rectifiability.  The main results of this paper concern the
properties of $N$-rectifiable sets, recovering many of the basic facts
about $\R^k$-rectifiable sets in $\R^n$ such as unique approximate tangent cones,
approximative qualities of tangent cones.  Precisely, we show the following
theorems.

\begin{MThm}  Let $N$ and $M$ be Carnot groups and suppose $E \subset
  M$ is a $N$-rectifiable set.  Then, 
\begin{itemize}
\item $E$ is $N$-approximable
\item For almost every $x \in E$, there exists a unique approximate
  tangent cone at $x$ which is isomorphic to $N$.
\end{itemize}
\end{MThm}

Roughly, a subset is $N$-approximable if at a.e. point there is a copy of $N$
sitting in $M$ which approximates $E$ locally in a measure theoretic
sense.  See definition 4.2 for a precise definition.  An approximate
tangent cone at $x \in E$ is a copy of $N=exp_{x_0}(V)$ (for some base
point $x_0$) in $M$ such that 
\[ \frac{\Ha^k_M(E \cap B_M(x,r)\minus X(x,V,s))}{r^k} \ra 0\] where $X(x,V,s)$
is roughly a ``cone'' over $N$ through a base point $x_0$ with ``slope''
$s$.  See section 5
for the precise definition.  Such notions are used
in \cite{Mattila} in the case of $\R^k$-rectifiability.  The reader
should note that this theorem is restricted to $N$ which are full
Carnot groups, not proper subgroups.  As mentioned above, this theorem, among other
applications, provides a class of Carnot-Carath\'eodory spaces which
are quite well behaved (with respect to local geometry).  The class of $N$-rectifiable submanifolds in a
given Carnot group $M$ has the property that each member is a
Carnot-Carath\'eodory manifold in its own right (by restriction of the
distribution) and that the tangent cone at almost every point is isomorphic to
$N$.  

Mirroring the Euclidean situation, we would like to show that the
rectifiability, approximability, and approximate tangent cone
conditions are equivalent.  Unfortunately, one cannot approach this
via the standard arguments because the various projections onto copies
of $N$ sitting in $M$ are almost never Lipschitz.  However, if we
define a strong approximate tangent cone at $x \in E$ to be a copy of
$N$ satisfying two conditions:
\begin{enumerate}
\item 
\[ \frac{\Ha^k_M(E \cap B_M(x,r) \minus X(x,V,s))}{r^{k \cdot depth(N)}} \ra 0 \] where
$depth(N)$ is the number of steps in the nilpotent Lie group $N$ and
$exp(V)=N$
\item For sufficiently small $s$ and $\epsilon>0$, \[X(x,V,s) \cap
X(x,V^\perp,s)\cap B_M(x,\epsilon) = \emptyset\]
\end{enumerate}

Then we have the following theorem.

\begin{MThm} Suppose $E$ is a $\Ha^k_M$-measureable subset of a Carnot
  group $M$ and $V$ is a subspace of $\mathfrak{m}$ so that $exp(V)$
 is isomorphic to $N$.  Further suppose that for $\Ha^k_M$ a.e. $y \in
 E$ $exp(V)$ is a strong
  approximate tangent cone at $y$.
Then, $E$ is $N$-rectifiable.
\end{MThm}

In \cite{Bell:1996}, Bella$\ddot{\i}$che poses the following
question:  if a Carnot-Carath\'eodory manifold has a locally uniform
tangent structure (i.e. if all tangent cones (in the sense of Gromov)
in a neighborhood are 
isomorphic as Lie groups), then does there exist a Lipschitz map
between the tangent cone at a point and a neighborhood of that point?
Our notion of rectifiability addresses this question
under the hypothesis that the Carnot-Carath\'eodory manifold is a
submanifold of a Carnot group.  Using theorem B, we see that if a
submanifold of a Carnot group has strong approximate tangent cones
isomorphic to $N$ almost everywhere (i.e. it has uniform tangent
structures in this measure theoretic sense), then there is a local Lipschitz
equivalence between $N$ and the Carnot-Carath\'eodory manifold.

Next, we turn to studying the rectifiability properties of
submanifolds in Carnot groups that do not necessarily enjoy tangent
structures which are Carnot groups themselves.  We focus on a class of
maps which appear in the work of Franchi, Serra Cassano and Serapioni
(\cite{FSSC2}), the $C^1_N$ maps.  There are functions $f:N \ra \R$
which have continuous horizontal derivatives (see the next section for
a precise definiton).  

\begin{MThm}  Let $N$ be a Carnot group of Hausdorff dimension $k$ and let $f:N \ra \R$ be a
  $C^1_N$ map.  Then, for a.e. $t \in \R$ and $\Ha^{k-1}_N$ a.e. $x
  \in Int(f^{-1}(t)  \minus \Sigma_t)$, there exists a subgroup
  $T_x \subset N$ and $r_x >0$ such that $B_N(x,r_x) \cap f^{-1}(t)$ is
  $T_x$-approximable at $x$ and $f^{-1}(t)$ has a unique approximate tangent
  cone isomorphic to $T_x$ at $x$.   
\end{MThm}

In the theorem, $\Sigma_t$ is the collection of characteristic points
of the level set $f^{-1}(t)$ and $Int(f^{-1}(t) \minus \Sigma_t)$
denotes the interior of $f^{-1}(t) \minus \Sigma_t$.  We note that
{\em a priori}, the interior may be empty. 
Note that in this theorem, we must allow $T$ to be a proper subgroup
of $N$ and that $T$, equipped with the metric from $N$ restricted to
$T$, may not be a Carnot group.  This theorem is a step
towards generalizing the techniques of geometric measure theory to the
Carnot setting.  In particular, this theorem gives information
concerning generalizing the notion of slicing manifolds by Lipschitz
maps.  Unfortunately, as shown by the limitations in this theorem (and
illustrated by an example in section 5) we cannot conclude that the
level sets are $N'$-rectifiable for some $N'$ - we lack a Lipschitz
map.  However, as evidenced by the theorem, many of the approximative
qualities of rectifiable sets are inherited by the level sets.  This
suggests a modification of the notion of rectifiable currents in the
Carnot setting based on these types of properties.  

The proofs of these theorems rest on extensions of Euclidean analytic
tools to the Carnot case.  The most useful one of these is the
(suitably defined) differentiability of Lipschitz maps on Carnot
groups which is originally due to Pansu (\cite{Pansu}).  The form of
the theorem used in this paper is an extension due to Vodopyanov and Ukhlov (\cite{VodUkh}) and recently proved using a different method by
Magnani (\cite{Magnani}).  In addition, in sections
2 and 3, we prove various lemmas concerning the properties of the
Hausdorff measure including a metric area formula.  Most of the proofs
of these lemmas are adaptations of arguments in Federer (\cite{Fed})
and the proof of the area formula follows Kirchheim's argument in
\cite{Kir}.  To prove theorems A and B, we follow arguments based on
arguments in \cite{Mattila} extended
using the lemmas and techniques described above.  However, to prove theorem C, we
diverge from the classical arguments, instead using smooth
approximations of the Lipschitz map (as in \cite{FSSC}, \cite{FSSC2}
and \cite{GarNh}) and apply metric arguments akin to those in
\cite{Pauls}.  

The structure of the paper is as follows:  section 2
reviews some of the known measure theory and differentiability results
for Carnot groups and proves many of the useful measure theoretic
lemmas needed in the proofs, such as properties of
Jacobians of maps and a weak Sard-like property.  Section 3 is devoted
to an area formula for Lipschitz maps between Carnot groups.  Sections
4-6 introduce $N$-rectifiability, $N$-approximability and prove
theorems A and B.  Section 7 is devoted to proving theorem C.  

\section{Background results for Carnot groups}
\subsection{Distances and measures}
Assume that $N$ is a connected, simply connected graded nilpotent Lie
group.  Recall that $N$ is graded if the Lie algebra decomposes as  $\mathfrak{n}=\mathcal{V} \oplus \mathcal{V}_2
\oplus ... \oplus \mathcal{V}_l$ where $[\mathcal{V}_1,\mathcal{V}_j]
= \mathcal{V}_{j+1}$.  
\begin{Def}We will call the integer $l$ the {\bf depth} of $N$
and denote it by \[depth(N)\]
\end{Def}  

We denote by $\mathcal{V}$ not only
the bottom level of the grading, but the left invariant vector bundle
generated by left translating $\mathcal{V}$ around $N$.  We also
assume that $\mathcal{V}$ is equipped with an inner product $\ip$,
which we also think of as a left invariant inner product on the
subbundle $\mathcal{V}$.  
\begin{Def}  A Carnot group is a quadruple $(N,\mathcal{V},\ip,d_N)$
  where $N$,$\mathcal{V}$ and $\ip$ are as above.  To define $d_N$ we
  let $\mathcal{H}$ be the absolutely continuous paths which are tangent almost everywhere to $\mathcal{V}$.  Then,
  the Carnot-Carath\'eodory distance is defined as
\[d_N(n_1,n_2) = \inf \left \{\int < \gamma', \gamma'>^\frac{1}{2}
  \bigg | \gamma
  \in \mathcal{H} \text{ and $\gamma$ connects $n_1$ to $n_2$} \right \} \]
\end{Def}

It follows from the construction that $d_N$ is a left invariant metric
on $N$ which admits a homothety, denoted $h_t$.  The homothety is
defined by its action on the Lie algebra where it acts on vectors in
$\mathcal{V}_i$ by multiplication by $t^i$.

We will often think of a Carnot group $N$ as identified with $\R^k$
via the exponential map and use an orthonormal basis $\{X_i\}$ for
$\mathfrak{n}$, adapted to the grading, as our coordinates for
$\R^k$. Given $n \in N$ and $r >0$, we define $Box(n,r)=\{exp_n((t_1,
\dots, t_k) | |t_j| < r^{d(j)} \}$ where $d(j)$ is the level of the
grading to which $X_j$ belongs.  Next, we recall the well
known Ball-Box theorem (see for example, \cite{Gromov:CC} or \cite{Mont}):

\begin{Thm}  Let $N$ be a Carnot group.  Then
  $\exists C >0$ such that, given $n \in  N$ and $r>0$,
\[ Box\left (b, \frac{r}{C}\right ) \subset B_N(n,r) \subset Box(n, Cr)\]
\end{Thm}

While the definition of $d_N$ is geometrically compelling, it is very
difficult to compute with.  Luckily, one can use any of a family of
quasi-norms on $N$ to aid in computation.  We describe one here.

\begin{Def} Suppose $N$ is a Carnot group with grading
  $\mathfrak{n}=\mathcal{V} \oplus \mathcal{V}_2 \oplus ... \oplus
  \mathcal{V}_l$.  For $n \in N$, decompose $n$ according to the
  Euclidean vector space basis for $\mathfrak{n}$,
  $n=e^{v_1 + v_2 + v_3 +... +v_l}$.  Then, we define a quasi-norm on
  $N$ as follows:
\[ |n|_{qn} = \left (\sum_{i=1}^{depth(N)} ||v_i||^\frac{2}{d(v_i)} \right
  )^\frac{1}{2} \] 
where $d(v_i)$ is the level of the grading of which $v_i$ is a
member and $|| \cdot ||$ is the Euclidean norm.  We also define a function on $N\cross N$ by $d_{qn}(n_1,n_2)
= |n_1^{-1}n_2|_{qn}$.
\end{Def}

Note that $d_{qn}$ is, by construction, left invariant and admits
$h_t$ as a homothety.  Thus, $d_N$ and $d_{qn}$ are biLipschitz
equivalent.  A good reference on quasi-norms on Lie groups is \cite{Goodman}.

To any metric on a metric space, we can introduce an associated
Hausdorff measure.  To fix notation, we define the Hausdorff measures
needed in this paper.

Given a Carnot groups $N$ with Carnot-Carath\'eodory metric $d_N$ and
a subset $U$ of $N$.  Let $\mathscr{E}_\delta$ be the collection of
closed $d_N$ balls of diameter less than or equal to $\delta$.  Then, 
\[\Ha_{N,\delta}^k(U) = \inf \{\sum diam_N(E_i)^k \mid  U \subset \cup
E_i, \; E_i \in \mathscr{E}_\delta \} \]
And, 
\[\Ha_{N}^k(U) = \lim_{\delta \ra 0} \Ha_{N,\delta}^k(U)\]

Note that we do not use a normalization by, for example, the measure
of the unit ball in $N$ and that we use covers by closed balls rather
than covers by any sets.  This Hausdorff measure is sometimes called
the spherical Hausdorff measure on $N$.  

\begin{Def} Let $N$ be a Carnot group.  The Hausdorff dimension of a
  set $E \subset N$ is
\[\sup \{s | \Ha^s_N(E)=0\}\]
\end{Def}

We quote a well known lemma (see, for example, \cite{Gromov:CC} 0.6B):
\begin{Lem}\label{hyper}
Let $S$ is a smooth hypersurface of a Carnot group $N$
where the Hausdorff dimension of $N$ is $k$.  Then, the Hausdorff
dimension of $S$ is $k-1$.  
\end{Lem}

We next define the notion of Ahlfors regularity.

\begin{Def}\label{Ahlfors}  Let $X$ be a complete metric space and $\mu$ a Borel
  measure on $X$.  We say $\mu$ is Ahlfors regular of dimension $s$ if
  there exists $C>0$ so that
\[C^{-1}r^s \le \mu(B_X(x,r)) \le C r^s\]
when $0 < r \le diam(X)$ and $x$ is a point where $\mu(B_X(x,r)) >0$
for all $r >0$.
\end{Def}

In this paper, we will be considering the Ahlfors regularity of
$\Ha_N^{k-1}$ when restricted to $S$, a hypersurface of $N$ (where $N$
satisfies the assumption of the most recent lemma).

Because Carnot-Carath\'eodory metrics and quasi-norms on Carnot groups are left
invariant, the Hausdorff measures associated to them are in fact all constant
multiples of the Haar measure on such a group. Thus, many of the same
measure theoretic results concerning densities, etc. are true for
$\Ha_{N}^k$ that are true for Lebesgue measure in $\R^m$.  For proofs
and discussions of these facts, in a more general setting, see
\cite{brokendreams}.    We will use the
biLipschitz equivalence of $d_N$ and $d_{qn}$ and the relation between
the Hausdorff measures freely in the computations below.

\begin{Lem} \label{lebdens} Let $k$ be the Hausdorff dimension of
  $N$.  If $U \subset  N$ is $\Ha_{N}^k$
  measurable then almost every point of $U$ is a Lebesgue density
  point.  In other words, for almost every $a \in U$,
  \[ \lim_{r \ra 0} \frac{\Ha^k_{N}(U \cap
  B_{N}(a,r))}{\Ha_{N}^k(B_{N}(a,r))} = 1\]
\end{Lem}

Next, we define the usual densities.

\begin{Def}
If $N$ is a Carnot group and $U\subset N$, $x \in N$, we define
\[ \Theta^{*,k}_N(U,x) = \overline{\lim}_{r \ra 0^+}\frac{\Ha^k_N(U
  \cap B_N(x,r))}{r^k} \]
and 
\[\Theta_{*,N}^k(U,x) = \underline{\lim}_{r \ra 0^+}\frac{\Ha^k_N(U
  \cap B_N(x,r))}{r^k} \]
If both of these coincide, then the common value is denoted $\Theta^k_N(U,x)$.
\end{Def}

Note that, due to the lack of normalization of the Carnot metrics and
the related Hausdorff measures, our estimates will be much less
precise than in the Euclidean case.

\begin{Lem} \label{densities}  If $U \subset  N$ and $\Ha_{N}^k(U) <
  \infty$ then there exists a constant $C>0$, depending only on the
  structure of the group $N$ and the Carnot-Carath\'eodory metric on
  $N$ so that
\begin{enumerate}
\item $C 2^{-k} \le \Theta^{*,k}_N(U,x) \le C$ for almost all $x \in U$.
\item If $U$ is $\Ha_{N}^k$ measurable then $\Theta^{*,k}_N=0$ for
  almost all $x \in N \minus U$.
\end{enumerate}
\end{Lem}

In this lemma (and in the rest of the paper), the phrase ``depends on
the structure of the group $N$'' means that the constant in question
depends only on the dimensions of the levels of the grading of
$\mathfrak{n}$.  
The proof of this theorem follows from the standard arguments in
$\R^n$ (see e.g. \cite{Fed}, 2.10.18 or \cite{Simon}, theorem 3.6)
combined with the observation that the Hausdorff measure 
derived from any Carnot-Carath\'eodory metric on a Carnot group is
left invariant and scales with the homothety.  

Next we quote two useful formulae.  Both can be found in \cite{CCcalc}
as equation 4.8 and proposition 4.9 respectively.

To state them, we assume $N$ is a Carnot group and pick a Riemannian
completion of the inner product defining the Carnot-Carath\'eodory
metric on $N$ which makes the grading orthogonal.

\begin{Def}  Let $N$ be a Carnot group and $\{X_1, \cdots, X_a\}$ be
  an orthonormal basis for $\mathcal{V}$.  
Let $\nabla f$ denote the Riemannian gradient of a
function $f$ and denote by $\nabla_0 f$ the horizontal gradient of
$f$,   
\[\nabla_0f = X_1 f \; X_1 + \cdots + X_af \; X_a \]
\end{Def}

\begin{Def} Let $S$ be a smooth surface in $N$, a Carnot group.
  Suppose further that $S$ is a level set of some smooth function
  $f$.  Then, $x\in S$ is called a {\bf characteristic point} if
  $\nabla_0f(x)=0$.
\end{Def}

Below $| \cdot |$ is the usual norm on $\R$.

\begin{Pro}[Coarea formula] \label{coarea} Suppose $f:N\ra \R$ is a smooth map from a
  Carnot group of Hausdorff dimension $k$ to $\R$ and $u$ is any
  nonnegative measurable function.  Then,
\[ \int_N u(x) |\nabla_0 f(x)| d\Ha^k_N(x) = \int_0^\infty
  \int_{f^{-1}(t)} u(y) d \Ha^{k-1}_N(y) dt \]
\end{Pro}

\begin{Pro} \label{surfmeas} Suppose $S$ is a level surface of a smooth real valued
  function $f$ on $N$, a Carnot group of Hausdorff dimension $k$.
  Then, for $\Ha^{k-1}_N$ almost everywhere on $S$,
\[d \Ha^{k-1}_N = \frac{|\nabla_0 f|}{|\nabla f|} dA\]
where $dA$ is the Riemannian area element on $S$.
\end{Pro}  

We note that the size of the characteristic locus has been
investigated by many authors recently.  In particular, Z. Balogh
(\cite{Zoltan}) has shown that one can construct a $C^{1,\alpha}$
surface in the Heisenberg group, $\mathbb{H}$, with a characteristic
locus of positive two-dimensional Euclidean Hausdorff  measure.  However, if the
surface is at least $C^1$ then the characteristic locus has
$\Ha^{k-1}_N$ measure
zero.  See, for example, \cite{Val} and \cite{Zoltan}.  

We end this section with a definition of the $C^1_N$ functions.  
\begin{Def}
If
$N$ is a Carnot group and $f:N \ra \R$ is a continuous function then
$f \in C^1_N$ if $\nabla_0 f$ exists (in the sense of distributions)
and is a continuous vector valued function.
\end{Def}

\subsection{Differentiability of Lipschitz maps}\label{diff}
We first quote a definition and a result of P. Pansu (see
\cite{Pansu}) crucial to our constructions:

\begin{Def}  A map $f:N \ra M$ between Carnot group is said to be
  differentiable in the sense of Pansu at $n \in N$ with differential
  $df_n$ if the limit 
\[ df_n(y) = \lim_{s \ra 0} h'_{\frac{1}{s}}f(n)^{-1}f(nh_s y) \] exists
and convergence is uniform for all $y \in N$.    
\end{Def}

The following theorem is a generalization of a theorem of Pansu
(\cite{Pansu}) due to Vodopyanov and Ukhlov (\cite{VodUkh}).
Recently, Magnani (\cite{Magnani}) has also proved this theorem using
a different technique.
\begin{Thm}  If $N$ and $M$ are Carnot groups and $f$ is a
  Lipschitz map from a measurable set  $U \subset N$ to $M$, then $f$ is
  differentiable almost 
  everywhere (in the sense of Pansu) and the differential $df_x$ is a
  graded group homomorphism at almost every point.
\end{Thm}

In \cite{Pauls}, the author proves a limited metric differentiability
for Lipschitz maps of Carnot groups into complete metric spaces.  In
the special case when the target space is another Carnot group, either
an extension of the author's arguments or an appeal to the result
above gives complete metric differentiability.  In some of the later
results, it is more convenient to work with the metric version of
differentiability, hence we state it here.

\begin{Thm}  Let $f: U \subset N \ra M$ be a Lipschitz map between Carnot
  groups.  Then $f$ is metrically differentiable almost everywhere.
  In other words, for almost every $n \in N$, 
\[ \Delta_n(y_1,y_2) = \lim_{t \ra 0} \left \{ \frac
  { d_M(f(n'h_ty_1),f(n'h_ty_2))}{t} \Bigg | n' \in B_N(n,t) \right \} \]
exists, the limit converges uniformly and $\Delta_n$ admits a
homothety and is left invariant under the action of $N$.  
\end{Thm}

The following lemma dictates exactly how the image of the differential
approximates the function at a point.  

\begin{Lem}\label{approx0}  Let $f: U \subset N \ra M$ be a Lipschitz
  mapping between Carnot groups, where $U$ is a positive measure subset of $N$ and let $df_x$ denote the Pansu
  differential at a point of differentiability, $x$.  Then, 
\[ d_M(f(xe^v),f(x)df_x(e^v)) = o(d_N(e^0,e^v)) \]
\end{Lem}
\pf This is a straightforward computation.  Let $v_0$ be the Lie algebra
vector in the direction of $v$ such that $d_N(e^0,e^{v_0})=1$ and let $t=d_N(e^0,e^v)$.
\begin{equation*}
\begin{split}
d_M(f(xe^v),f(x)df_x(e^v)) &=
d_M(f(x)^{-1}f(xh_te^v),df_x(h_te^{v_0})) \\
&= d_M(h_tdf_x(e^{v_0}),df_x(h_te^{v_0})) + o(t) \\
&= o(t)
\end{split}
\end{equation*}
The last equality stems from the fact that Pansu's differential
intertwines the homotheties of $N$ and $M$.  $\qed$
\newline \newline

\noindent
{\em Remark: }  One should note that the resulting metric $\Delta_x$
is very close to being a Carnot-Carath\'eodory metric on $N$ (the only
possible degeneracy is that
the inner product on each fiber may be only semi-definite).  To see
this, one combines the fact that the differential is a graded group
homomorphism and the definition of the metric $\Delta_x$ with the
previous lemma, concluding that one may recognize $\Delta_x(n_1,n_2)$
by the quantity $d_M(df_x(n_1),df_x(n_2))$.   As we shall see in
the next section, at points of differentiability when the map is
suitably nondegenerate, $\Delta_x$ is a well defined
Carnot-Carath\'eodory metric.

\subsection{Jacobians and a weak Sard-like theorem}
Let $f:U \subset N \ra M$ be a Lipschitz map between a positive
measure subset of a Carnot group and another Carnot group and let $k$
be the Hausdorff dimension of $N$.  

\begin{Def}\label{Jac}  For each $x \in U \subset N$ we define the Jacobian of
  the map $f$ at $x$ by
\[ J(x) = \lim_{t \ra 0} \left  
\{ \frac{\Ha^k_M(f(B_N(y,t)))}{\Ha^k_N(B_N(y,t))} \Bigg | y \in
  B_N(x,t) \right \} \] 
\end{Def}

By arguments analogous to those in \cite{Pauls} concerning the
existence of the metric differential, $J(x)$ exists almost everywhere
and scales appropriately.
Moreover, we now prove that the image of the set where
$J(x)=0$ has $\Ha^k_M$-measure zero.

\begin{Lem} Let $Z=\{x \in U \subset N | J(x)=0\}$.  Then, $\Ha^k_M(f(Z))=0$.  
\end{Lem}
\pf Fix $\epsilon >0$, $n_0 \in N$ and $0 < R < \infty$.  Let $Z_{R}=Z
\cap B_N(n_0,R)$.  We will show that $f(Z_{R})$ has 
measure zero.  We may cover
$Z_R$ by balls $B_N(n,r)$ with the property that $n \in Z_R$ and
$\frac{\Ha^k_M(f(B_N(n,r)))}{\Ha^k_N(B_N(n,r))} < \epsilon$.
By the assumption that $J(n)=0$, this covering is fine, so using the
Vitali covering lemma, we refine the cover to a countable disjoint
collection of balls $\{B_i=B_N(n_i,r_i)\}$ with $r_i < R$ which cover
almost all of 
$Z_R$.  Since $f$ is Lipschitz, we know that $\cup_i f(B_i)$ covers
$\Ha^k_M$ almost all of $f(Z_{R})$ as well.  Thus 
\begin{equation*}
\begin{split}
\Ha^k_M(f(Z_{R})) &\le \sum_i \Ha^k_M(f(B_N(n_i,r_i))) \\
&< \epsilon \sum_i \Ha^k_N(B_N(n_i,r_i)) \\
&\le \epsilon \Ha^k_N(Z_R) \\
&\le \epsilon \Ha^k_N(B_N(n_0,R))
\end{split}
\end{equation*}

Thus, since $\epsilon$ is arbitrary, $\Ha^k_M(f(Z_{R})) =0$.  The
result now follows easily.  $\qed$

Next we prove a lemma analogous to Sard's theorem concerning the
measure of the image of the ``degenerate'' set.  Suppose $f:U \subset N \ra M$ is a
Lipschitz map and let $N_x = f(x)df_x(N)$.  
\begin{Lem}\label{Sard}  Suppose $\dim_\Ha(N) \le \dim_\Ha(M)$.  Then,
  $\Ha_M^k(\{f(x) | \Ha_M^k(N_x) = 0\})=0$. 
\end{Lem} 
\pf To begin with, we may assume that at all points Pansu's derivative exists and $J(x)$ exists for all $f(x)$ we are considering since the set of the
complementary points has measure zero.  First we observe that using
the definition of Pansu's differential, the uniformity of its
convergence, the left 
invariance and homothety of the Carnot-Carath\'eodory metric $d_M$, we
have that $d_M(f(xh_te^v),f(x)h_tdf_x(e^v)) =
o(t)$.  In particular, this says that $f(B_N(x,t))$ lies in a $o(t)$
neighborhood of $f(x)df_x(B_N(e^0,t))$.  Since $f(B)$ has Hausdorff
dimension less than or equal to $k$, we know that
$\Ha^k_M(f(B_N(x,t))) \le \Ha^k_M(df_x(B_N(x,t))) 
+ o(t^k)$.  Now consider an image point
$f(x)$ such that $\Ha_M^k(N_x) = 0$.  We discuss separately two cases:
ker$(df_x)=\{e^0\}$ and ker$(df_x) \neq \{e^0\}$.  In the first case,
this together with the fact that $df_x$ is a graded group homomorphism
implies that the image of $df_x$ is isomorphic to $N$.  Thus,
$\Delta_x$ is positive definite on $N$ and, since $\Delta_x$ is left
invariant and admits a homothety, $df_x: (N,\Delta_x) \ra (M,d_M)$ is
biLipschitz onto its image and hence cannot have $\Ha^k_M$-measure
zero.  So, we may assume that the kernel is nontrivial.  Hence, $N_x$
is isomorphic to a quotient of $N$ (by the kernel), call it $N'$.  In
particular, 
realizing $N'$ as a subgroup of $N$, we see that $\Ha^k_N(N') = 0$
and so, since $df_x$ is Lipschitz, $\Ha^k_M(N_x) \le L^k
\Ha^k_N(N')=0$.  Thus, we have that $\Ha^k_N(f(B_n(x,t))) = o(t^k)$,
yielding the desired result.  $\qed$

%   Now, consider the map $F:N \ra M \cross N$ given by
% $F(n)=(f(n),h_\epsilon n)$.  One quickly calculates that
% \[ \epsilon d_N(n_1,n_2) \le d_{M\cross N}(F(n_1),F(n_2)) \le
% \sqrt{L^2 + \epsilon^2}\epsilon d_N(n_1,n_2) \]
% At a point $x$ where $df_x$ has a kernel, we compute the Jacobian of
% $F$ at $x$ by estimating the $k$ dimension Hausdorff measure of
% $F(B_N(x,t))$ with respect to the standard quasi-norm, finding that
% $\Ha^k_{qn}(F(B_N(x,t))) \le C \epsilon t^{k-1}$.  Since the
% projection onto the first factor of $M \cross N$ is $1-$Lipschitz, we
% see that $J(x) \le  C \epsilon t^{k-1}$ as well.  Thus, since
% $\epsilon$ is arbitrary, we have that $\{ f(x) | \Ha^k(N_x)=0 \} = \{ f(x) | J(x)=0\}$ and so, by lemma \ref{Jac}, we have the desired
%   conclusion.  $\qed$.

\section{An area formula}
In this section, we provide a proof of an area formula for
Lipschitz maps.  A change of variables formula is proved by Vodopyanov and Uhklov
(\cite{VodUkh}) using methods based on Pansu's techniques in
\cite{Pansu}.  One should be able to extend their arguments to prove
the statement below.  Also, recently, Magnani (\cite{Magnani})
independently proved the same area formula using a different (but
equivalent) definition of Jacobian.  We include the proof here for completeness.  Again, let $f:U \subset N \ra M$
be a Lipschitz map of a Carnot group $N$ to a Carnot group $M$.  Let
$k$ be the Hausdorff dimension of $N$.  

\begin{Lem} \label{Jac approx}  Fix $\lambda >1$ and let $E$ be a
  measurable set in $N$ such that for every $x \in E$, $\Delta_x$
  exists and is nondegenerate.  Then there exists a countable Borel cover of
  $E$, $\{B_i\}$, and for each $B_i$, a left invariant CC metric,
  $d_i$, on $N$ such that for all $x,y \in B_i$,
\[ \lambda^{-1} d_i(x,y) \le d_M(f(x),f(y)) \le \lambda d_i(x,y)\]  Moreover, 
\[\lambda^{-k}
  \frac{\Ha^k_i(B_N(x,1))}{\Ha^k_N(B_N(x,1))} \le J(x)\le \lambda^{k}
  \frac{\Ha^k_i(B_N(x,1))}{\Ha^k_N(B_N(x,1))}\] for all $x \in B_i$
  which are density points of $B_i$ with respect to the Hausdorff
  measure, $\Ha_i^k$, constructed with respect to the distance $d_i$.
\end{Lem}

\pf This theorem follows exactly as lemma 4 in \cite{Kir} which,
in turn, follows much of the argument of lemma 3.2.2 in \cite{Fed}.
The metrics $d_i$ are given by $d_i(n_1,n_2)=d_M(df_x(n_1),df_x(n_2))$.

%\begin{Cor} For a.e. $x \in E$, $J(x)=J_{\Delta_x}(x)$.
%\end{Cor}

While we explore the concept of $N$-rectifiability in detail in the
next section, we define a set $E \subset M$ to be $N$-rectifiable if
it is the Lipschitz image of a positive measure subset of $N$.  

\begin{Cor}\label{denseone} If $E \subset M$ is $N$-rectifiable, then
  for almost every $x \in E$, $\Theta^k(\Ha^k_M \hook E,x)=\Ha^k_N(B_N(e^0,1))$.
\end{Cor}

\pf Using lemma \ref{Jac approx} above, this follows at density points
of the $B_i$.  More precisely, fix $\lambda >1$ and let $B_i$ and
$d_i$ be as in the lemma.  Now, for each point of density of $B_i$
which is also a point of Pansu differentiability, let $K_i$ be the
preimage of $B_i$.  We have,  
\begin{equation*}
\begin{split}
\lambda^{-2k} \frac{\Ha^k_N(K_i \cap
  B_N(x,\frac{\delta}{\lambda}))}{\left ( \frac{\delta}{\lambda}
  \right )^k}
  &\le \frac{\Ha^k_M(B_i\cap B_M(f(x),\delta))}{\delta^k} \\
&\le \lambda^{2k} \frac{\Ha^k_N(K_i \cap
  B_N(x,\delta\lambda))}{( \delta\lambda )^k}
\end{split}
\end{equation*}

Thus, taking advantage of the assumption that $x$ is a density point
and letting $\delta$ go to zero and $\lambda$ go to one, we have the
density is constant almost everywhere.  The result follows.  $\qed$

We note that, if we had defined $\Ha_N^k$ with a normalizing constant,
we would then have the analogue of the classical fact that an
$N$-rectifiable set has density one almost everywhere.  

\begin{Thm}\label{area}  Suppose $f:N \ra M$ is a Lipschitz map.
  Then, for any $\Ha^k_N$-measurable set $E$,
\[ \int_E J(n) d \Ha^k_N(n) = \int_M N(f|_E,m) d \Ha^k_M(m) \]
\end{Thm}

\pf Assume for a moment that at every point in $E$, the Pansu
differential is well defined and has trivial kernel. In other words, $E \subset \{x | \Delta_x \text{ is a nondegenerate CC metric}\}$.  Then, fixing
$\lambda >1$ and using lemma \ref{Jac approx}, we find a countable
cover $\{B_i\}$ with the approximative properties described in the
lemma.  Let $J_i=\frac{\Ha^k_i(B_N(x,1))}{\Ha^k_N(B_N(x,1))}$.
Note that, under the assumption of nondegeneracy, $\Ha^k_N$ and
$\Ha^k_i$ are both constant multiples of one another and $\Ha^k_i=J_i
\Ha^k_N$.  Using these facts, we have that  
\begin{equation*}
\begin{split}
\lambda^{-k} J_i &\le J(x)\le \lambda^{k} J_i \text{ for all $x
  \in B_i$}\\
\intertext{and hence,}
\lambda^{-2k} \int_{B_i \cap E} J(n) d\Ha^k_N(n) &\le \lambda^{-k} \int_{B_i \cap E}
  J_i d\Ha^k_N(n) = \lambda^{-k} \int_{B_i \cap E} d\Ha^k_i\\& \le
  \int_M N(f|_{B_i \cap E},m) d \Ha^k_M(m) \le \lambda^k \int_{B_i\cap
  E} d\Ha^k_i \\ &=
  \lambda^{k} \int_{B_i \cap E}  J_i d\Ha^k_N(n) \le \lambda^{2k}
  \int_{B_i \cap E} J(n) d\Ha^k_N(n)\\
\end{split}
\end{equation*}

Thus, summing over $i$ and letting $\lambda \ra 1$ we have the desired
result.  

By lemma \ref{Sard} and the arguments used to prove it, we see that at
points $x$ of $E$ where $\Delta_x$ is not defined or is degenerate,
$J(x)=0$ and the set of all such $x$ maps to a $\Ha^k_M$
measure zero set, making both sides of the desired equation zero.
$\qed$.

In light of the discussion above, we note that we could have defined
the Jacobian via the Pansu differential as follows:  at a point $x$ of
Pansu differentiability, $J(x)=
\frac{\Ha^k_M(df_x(B_N(e^0,1)))}{\Ha^k_N(B_N(e^0,1))}$.

\section{CC-rectifiability}
Next, we introduce the definitions and basic properties of a theory of
rectifiability for subsets of Carnot groups.  One
should view this analogously to the Euclidean case:  Euclidean
rectifiable sets are viewed as sets which have, in some sense, a
manifold structure while CC-rectifiable sets will have a
``manifold'' structure where the local geometry
modeled by general Carnot groups rather than simply by $\R^n$.  

We begin with the relevant definitions.  Assume that $N$ and $M$ are
Carnot groups and that the Hausdorff dimension of $N$ is $k$.

\begin{Def}\label{CCrect}  Let $N$ be a Carnot group.  A subset $E$ of
  another Carnot Group  $(M,d_M)$ is said to be {\bf $N$-rectifiable} if there exists
  $U$ an positive measure subset of $N$ and a Lipschitz map $f:U \ra
  M$ such that 
  $\Ha^k_M(E \minus f(U))=0$.  $E$ is said to be {\bf countably
  $N$-rectifiable} if there exist a countable number of $U_i$ and
  $f_i:U_i \ra Y$ Lipschitz  with $\Ha^k_M(E \minus \cup_i
f_i(U_i))=0$.
\end{Def}

As with the case of $\R^n$-rectifiability, we will be developing the
notion of approximate tangent cones and their
relation to rectifiability.  
Next, we wish to consider when a subset of $M$ is well approximated by
$N$.  If $i:N \ra M$ is a graded group homomorphism, let
$N(t)= \{ x | d_M(x,i(N)) \le t\}$.  

\begin{Def}\label{Napp}  A subset $E$ of $M$ is {\bf $N$-approximable}
  at a point $a \in E$  if, for 
  $\alpha >0$, there exist an injective homomorphism $i:N \ra M$
  with $i(e^0)=e^0$ and constants $r_0 >0, \theta >0$ such that for any
  $0 < r< r_0$,
\begin{equation}\label{approx1}
\Ha^k_M(E \cap B_M(b,\alpha r)) \ge \theta r^k \text{for $b \in a
  \cdot i(N)\cap
  B_M(a,r)$} 
\end{equation}
and
\begin{equation}\label{approx2}
\Ha^k_M(E \cap B_M(a, r) \minus a \cdot N(\alpha r)) < \alpha r^k
\end{equation}

A set $E$ is called {\bf $N$-approximable} if $\Ha^k_M$ a.e. point in
$E$ is $N$-approximable.
\end{Def}

Next, we have the first theorem, which follows from the
differentiability of Lipschitz maps described in theorem \ref{diff}.
The proof of the theorem follows the argument of theorem 15.11 in
\cite{Mattila} but with the appropriate changes for the Carnot case.
The reader should also note that another key component of the proof
relies on the fact that a Carnot group and its tangent cone at a point may be
identified due to the existence of a homothety of the
Carnot-Carath\'eodory metric thus allowing
the homomorphism of tangent cones to be translated as a statement
concerning sets in the group.  

\begin{Thm}\label{rect/approx}  Every countably $N$-rectifiable $E \subset M$
  with nonzero $k$-dimensional Hausdorff measure is $N$-approximable.
\end{Thm}
\pf In this proof, the main information to keep in mind is that the differentiability
theory of Lipschitz maps between Carnot groups is close enough to that
of Lipschitz maps between Euclidean spaces, allowing many of the
Euclidean arguments to be used.  

Fix $\epsilon >0$.  First, we reduce to an easier case.  Let $0 < \alpha <1$.  Since $E$ is
countably $N$-rectifiable, we consider an L-Lipschitz map between $E'
\subset N$ such that $f(E') \subset E$.  To make the reduction to
nicer sets, where the lower density is bounded away from zero, we use
lemma \ref{denseone}.

Using this, we can cover $E'$ up to a set of measure zero by a
countable union of subsets, $S$, with the following properties.  For
each $S$, there exists constants $\theta >0$, $r_0 >0$ such that for
$x \in S$ and $0< r < r_0$, 
\begin{equation} \label{dens1}
\Ha^k_M(f(E') \cap B_M(x,r)) \ge \theta r^k
\end{equation}  On these sets, we shall verify the properties of
$N$-approximability.  

To create a approximating isometric embedded copy of $N$, we use
Pansu's differential mapping.  Given a point $x$ at which $f$ is
differentiable, Let $i_x(y)=f(x)df_x(x^{-1}y)$ and $N_x=i_x(N)$.  Note
that the multiplication 
is multiplication in the Carnot group $M$.  Since $df_x$ is a group
homomorphism, we know that $N_x$ is a subgroup of $M$ and from lemma \ref{Sard}
we know that for almost every such $x$, $N_x$ has large
$k$-dimensional Hausdorff measure.  In
other words, the mapping $df_x$ has, in some sense, full rank.
Interpreting this result, we now show that at almost every image point 
of $f$, $N_x$ is a isomorphic copy of $N$.  If it were not isomorphic, it would
be isomorphic to a subgroup $N'$ of $N$ (since $df_x$ is a graded
homomorphism). 
As a graded subgroup of $M$, $N_x$ inherits a Carnot-Carath\'eodory metric
which is biLipschitz to the restricted Carnot-Carath\'eodory metric on $N'$
because all Carnot-Carath\'eodory metrics given by varying the norm on
the distribution of a single nilpotent
Lie group are biLipschitz to one another.  Thus, $N'$ and $N_x$ have
strictly lower Hausdorff dimension than $N$ and so lemma \ref{Sard}
shows this may only happen at almost every $f(x)$.  Once again, we
replace the set $E'$ with a full measure subset of points $x$ such
that $N_x$ has ``full rank'' in the sense described above.
Thus, for each $x$, $(N_x, d_M|_{N_x})$ is biLipschitz to $(N,d_N)$.
Denote by $l(x)$ the lower Lipschitz constant for each $x$.  By
construction, $l(x)>0$ for all $x \in E'$.  

Recalling that the Lebesgue density
theorem holds in $N$ for the measure $\Ha^k_N$ and using the
approximations developed above, we will now find numbers $r_0>0$ and
$\delta >0$ and a compact set $E_0
\subset E'$ with $\Ha^k_N(E' \minus E_0) < \epsilon$ consisting only
of density points of $E'$ with nice approximative properties.  We
begin by picking a $\delta < \min\{\frac{\alpha}{4},\frac{1}{L}\}$.  First,
pick a compact subset of the density points of $E'$, $E_0$, so that for $x \in E_0$, $l(x)
\ge 2 \delta$ ({\em Property $1$}).  Second,  because $N_x$ approximates $f(E')$ well at
$f(x)$ by lemma \ref{approx0}, I can pick $r_0$ small enough so that for $x \in E_0$,
\[d_M(f(y),i_x(y)) < \delta^2 d_N(x,y) \text{ for $y \in B_N(x,r_0)$
  \hfill \; \; \; ({\em Property $2$})}\]
So far, the choice of $r_0$ and $E_0$ depends on $\delta$ only.  Last,
we pick $\delta$ and, possibly readjusting $r_0$ to be smaller still,
we may guarantee that for 
$x \in E_0$, $0< r < r_0$ and $y \in B_N(x, \frac{r}{\delta})$, 
$d_N(y,E_0) < \delta^2 r$ ({\em Property $3$}).  Property $3$ follows from the fact that all points
in $E_0$ are density points of $E'$.   

Now, writing $E_0$ as the union of finitely many subsets $E_i$ with
$diam_N(C_i)< r_0$, we examine each $C_i$ individually.  Consider a
point $f(x)$ with $x \in C_i$ and $\Theta^k_M(f(E') \minus
f(C_i),f(x))=0$.  Since almost every point in $f(C_i)$ has this
property, we will consider only these points.  Let $0 < r <
\frac{\delta r_0}{2}$ and pick $i_x(y) \in N_x \cap B_M(f(x),r)$.
Note that $ y \in B_N(x,\frac{r}{\delta})$ by the property  $1$.  By
property $3$ above, there exists $z \in E'$
such that $d_N(y,z) < \frac{2 r}{\delta}$.  Using the triangle
inequality, the fact that $i_x$ is $L$-Lipschitz (recall that $L <
\frac{1}{\delta}$)  and property $2$, we get 
\[d_M(f(z),i_x(y)) \le \delta^2 d_N(x,z) + L d_N(y,z) \le 3r\delta \]

Since $4 \delta < \theta$, we have from equation \ref{dens1} that
\[ \Ha^k_M(f(E') \cap B_M(i_x(y),\alpha r)) \ge \Ha^k_M(f(E') \cap
B_M(f(z),\delta r)) \ge \theta \delta^k r^k\]

Taking $\lambda = \theta \delta^k$ we have verified equation
\ref{approx1}.

To verify equation \ref{approx2}, we observe the following:

First, by property $2$,
 \[f\left(C_i \cap B_N\left (x, \frac{r}{\delta}\right)\right) \subset N_x(\delta r)
\subset N_x(\alpha r)\]
Second, using the lower Lipschitz bound on $C_i$ and 
property $2$ again, 
\[ f\left(C_i \minus B_N\left(x, \frac{r}{\delta}\right)\right) \subset M \minus
B_M(f(x),r)\]
and so, $f(C_i) \cap B_M(f(x),r)$ lies inside $N_x(\alpha r)$.  Given
the density assumption on $f(x)$, that $\Theta^k_M(f(E') \minus
f(C_i),f(x)) = 0$, this implies that, again possibly shrinking $r_0$,
that equation \ref{approx2} holds.  $\qed$

\section{Implications of Rectifiability}

To do this, we need to use a type of
projection mapping analogous to the projections onto planes used in
the usual theory.  

\begin{Def}\label{proj}  Let $V$ be a vector subspace of
  $\mathfrak{n}$ and let $<\cdot ,\cdot >$ be a Riemannian completion
  of the CC-inner product on $\mathfrak{n}$ which makes the grading
  orthogonal.  Denote by $v^\perp$ the orthogonal complement of $V$ (with
  respect to this inner product).  Let $\mathfrak{pr}_V:
  \mathfrak{n} \ra V$ by the projection of $\mathfrak{n}$ onto $V$ and
  let $P_V:N \ra e^{V}$ be the map $exp \comp \mathfrak{pr}_V \comp
  exp^{-1}$ where $exp$ is the usual exponential map.  Also, let $Q_V$
  be the map $P_{V^\perp}$.  In each of the applications of the
  projection mappings below, there is an understood base point for the
  exponential map.  
\end{Def}

It is a direct consequence of the equivalence of the
Carnot-Carath\'eodory metric and $d_{qn}$ that the projection defined above is a Lipschitz
map only very rarely.  Although the
map may not be Lipschitz, the map defined 
above is a projection in the sense that $P_V\comp P_V$ is the identity map.  
One can easily construct examples of projections which are not
Lipschitz.  Note, however, that these projections are
H\"older with exponent $\frac{1}{depth(N)}$, where
$depth(N)$ is the numbers of levels in the grading of $\mathfrak{n}$,
when we restrict our view to a neighborhood of the origin.
\newline \newline

\noindent
{\em Example: } Consider the 3-dimensional Heisenberg
  group, $H^3$, with a left invariant Carnot-Carath\'eodory metric, $\cc$.
  Let $V$ be the vector subspace spanned by the ``Z''
  (non-distributional) direction.  Then, the exponential image of $V$
  in $H^3$ is a 1-parameter subgroup but $V$ is not a Lie subalgebra
  with compatible grading.  We will now directly show that this
  projection is not Lipschitz.  Not surprisingly, this stems from the
  nontrivial bracket structure.  Giving $H^3$ coordinates $\{X, Y,
  Z\}$, we describe an element $e^{aX+bY+cZ}$ by the triple
  $(a,b,c)$.  Then $\cc((a,b,c),(\alpha,\beta,\gamma)) =
  \cc((0,0,0), (\alpha-a,\beta-b,\gamma-c +\frac{1}{2}(\alpha b -
  a \beta)))$.  Fixing $\epsilon >0$ and taking $\alpha = a +
  \epsilon$, $\beta = b + \epsilon$ and $a=\frac{2(\gamma-c) +b
    \epsilon}{\epsilon }$, this simplifies to $\cc((0,0,0),
  (\epsilon,\epsilon, 0))= \epsilon \cc((0,0,0),(1,1,0))$.  Computing
  the distance under the image of 
  $P_V$, we have \[ \cc(P_V((a,b,c)), P_V((\alpha,\beta,\gamma))) =
  \cc((0,0,c),(0,0,\gamma)) = \sqrt{|\gamma -c|}\cc((0,0,0),(0,0,1))\]
  Thus, since $\epsilon, \gamma$ and $c$ are arbitrary, we see that
  the map cannot be Lipschitz.  Notice, however, that if the two
  points are always within some fixed ball around the origin, this map
  is H\"older with exponent $\frac{1}{2}$ $\qed$

 \begin{Def}\label{cone}  Let $V$ be a subspace of
   $\mathfrak{m}$.  If $m_0 \in M$, $0< s <1$, and $0 < R< \infty$,
   we define \[ X(m_0,V,s)= \{ m_1 \in M | d_M(Q_V(m_1), Q_V(m_0)) < s
   d_M(m_1,m_0) \} \]  Moreover, let \[X(m_0,R,V,s)= X(m_0,V,s) \cap
   B_M(m_0,R)\]
 \end{Def}

 The reader should look closely at this definitions of ``cones''.
 While they are the same in form as the Euclidean versions used by
 Mattila in \cite{Mattila}, the actual objects look somewhat different
 in practice (it is useful and relatively easy to investigate the shape
 of these objects in the Heisenberg group).  While they still have the
 same general form of an ``X'' 
 emanating from the base point, the spreading of the limbs of the ``X''
 is no longer quite linear.  In practice, we will see that this will
 make no difference for our applications.

Next, we assume that $N$ and $M$ are Carnot
groups of Hausdorff dimensions $k$ and $l$ respectively.
\begin{Def} Suppose $Y \subset M$, $m \in M$ and $V$ is a subspace of
  $\mathfrak{m}$ with $exp(V)$ isomorphic to $N$.  We say that $exp_m(V)$ is
  an approximate tangent cone for 
  $Y$ at $m$ if $\Theta^{*k}(Y,m) > 0$ and for all $0 < s <1$,
\[\lim_{r \ra 0^+} \frac{\Ha^k_M(Y \cap B_M(m,r) \minus
  X(m,V,s))}{r^k} = 0\]

Following the Euclidean notation, we denote the set of all approximate
tangent cones as $apTan^N(Y,m)$.  Note that this set depends on a
choice of Carnot group $N$.  
\end{Def}

Next, we come to the first main theorem; we show that
$N$-rectifiability implies $N$-approximability and the existence
of approximate tangent cones almost everywhere.

 \begin{Thm} \label{Equiv}  Let $N$ and $M$ be Carnot groups and the Hausdorff
   dimensions of $N$ is $k$.  Let $Y$ be an $N$-rectifiable subset of $M$ with $\Ha^k_M(Y) < \infty$.
   Then,
\begin{enumerate}
%\item $Y$ is $N$-rectifiable 
 \item $Y$ is $N$-approximable.
 \item For $\Ha^k_M$ almost all $y \in Y$, there is a unique
   approximate tangent cone at $y$ isometric to $N$.
% \item For $\Ha^k_M$ almost all $y \in Y$, there exists some  
%   approximate tangent cone at $y$ isometric to $N$.
 \end{enumerate}
 \end{Thm}

\pf As with the previous results, this proof follows the development
of \cite{Mattila} with changes for the Carnot case.  

\noindent
$(1)$ is theorem 4, so it remains to prove the other implication.

\noindent
{\bf $(1) \Rightarrow (2)$} \newline
Notice that, for $\epsilon >0$ sufficiently small, $B_M(y,r)
\minus X(y,N,s) \subset (B_M(y,r) \minus N(\epsilon s r)) \cup
B_M(y, \epsilon r)$.  From this, and the assumption that $Y$ is
$N$-approximable (use the second property in the definition), we see that a copy of $N$ is an approximate tangent cone for
$Y$ at almost every point.  To see that the approximate tangent cone
is unique, we use property $(1)$ in definition 4.2 and the definition
of approximate tangent cone. $\qed$

\section{Measures in cones} 
Theorem \ref{Equiv} stops short of proving the equivalence of
$N$-rectifiability, $N$-approximability and the existence of
approximate tangent cones.  However, we will show that
$N$-rectifiability is implied by the existence of a stronger notion of
approximate tangent space.  The notion of
approximate tangent cone in the last section is not strong enough to
easily characterize 
rectifiable sets due to the failure of the projections to be
Lipschitz.  In this section, we introduce a stronger version of
approximate tangent cones. 

\begin{Def} Let $M$ be a Carnot group.  Suppose $Y \subset M$, $m \in M$ and $V$ is a subspace of
  $\mathfrak{m}$ with $exp(V)$ isomorphic to $N$, where $N$ is a
  subgroup of $M$.  We say that
  $exp_{m_0}(V)$ (for some $m_0$ in $M$) is a strong approximate
  tangent cone for 
  $Y$ at $m$ if $\Theta^{*k}(Y,m) > 0$ and for all $0 < s <1$,
\[\lim_{r \ra 0^+} \frac{\Ha^k_M(Y \cap B_M(m,r) \minus
  X(m,V,s))}{r^{k depth(N)}} = 0\]
and that for sufficiently small $s$ and $\epsilon>0$, \[X(m,V,s) \cap
X(m,V^\perp,s)\cap B_M(m,\epsilon) = \emptyset\]

We denote the set of all strong approximate tangent cones as $SapTan^N(Y,m)$.
\end{Def}

\noindent
{\em Remarks:} \begin{enumerate}
\item  This definition is highly restrictive and is certainly
not the sharpest definition to have the consequence of
$N$-rectifiability. 
\item Note that the group $N$ in the previous definition is now not
  restricted to being a Carnot group itself, it is merely a subgroup
  of a Carnot group.
\end{enumerate}

We will need a
method of producing Lipschitz maps.  To do so, we follow the idea used
in \cite{Mattila} of considering measures of the
intersections of cones and the set $E$ to conclude $N$-rectifiability.

 We begin with a simple lemma to aid in producing Lipschitz maps to
 rectifiable sets.  

 \begin{Lem}  Suppose $Y \subset M$ and $V$ is the orthogonal
   complement of a subspace $\mathfrak{n}$ of $\mathfrak{m}$,
   $0< s< 1$ and $0 <  R<\infty$ then if $Y \cap X(y,R,V,s)= \emptyset$
   for all $y \in Y$,  then $Y$ is $N$-rectifiable.
 \end{Lem}

 \pf Geometrically, this says roughly that the Lie algebra preimage of
 the set $Y$ at each point locally lies in $\mathfrak{n}$.  Thus it
 makes sense that one can approximate the set by Lipschitz images of
 pieces of $N$.  To prove this, we first may assume that $diam_M(Y)<
 R$, otherwise we may simply cut up $Y$ into a countable number of such
 pieces.  We note that, for this lemma to be true, we must assume that
 the maps $Q_V$ are all the same mapping, i.e. that the projections
 are defined via the exponential map at the {\em same} basepoint (it
 does not matter now what basepoint we pick).
 Next, we observe that if $y_1,y_2 \in Y$, the hypothesis of 
 empty intersection implies that $d_M(Q_V(y_1),Q_V(y_2)) > s
 d_M(y_1,y_2)$.  Thus, the map $(Q_V|_Y)^{-1}$ is
 Lipschitz with constant less than $\frac{1}{s}$.  We see that there exists
 $\overline{Y} \subset N$, a positive measure subset such that
 $(Q_V|_Y)^{-1}: \overline{Y} \ra Y$.  Hence, $Y$ is
 $N$-rectifiable. $\qed$

 The reader should note that one implication of the lemma is that a
 purely $N$-unrectifiable set must have a nonempty intersection with
 the cone for
 almost all points in the set.  

Next, we refine the lemma above, replacing the requirement of empty
intersection with a bound on the measure of the intersection.

 \begin{Lem} \label{meascone}  Let $V, \mathfrak{n}$ and $s$ be as
   above.  Let $0 < \delta < \infty$ and $0 < \lambda < \infty$.  Let
   $k$ be the Hausdorff dimension of $N$.  If $Y \subset M$ is purely
   $N$-unrectifiable and $\Ha^k_M(Y \cap X(y,r,V,s)) \le \lambda (r
   s)^{k depth(N)}$ for $ y \in Y$ and $0 < r< \delta$, then $\Ha^k_M(Y \cap
   B_M(m,\frac{\delta}{6})) \le C \lambda \delta^k$ for all $m \in M$
   where $C$ is a constant depending only on $k$ and the structure of
   the Lie group $M$.
 \end{Lem}

 \pf  This proof is similar to the proof of lemma 15.14 in
 \cite{Mattila}.  Through the usual reduction of $Y$ via decomposition
 into countably many sets, we may assume that $Y \subset
 B_M(m,\frac{\delta}{6})$.  Recall also that $Q_V$ is a locally 
 H\"older map with exponent at worst $\frac{1}{depth{M}}$.  Let
 $c_0$ be the H$\ddot{\text{o}}$lder constant.  By the previous lemma,
 we know that $Y \cap X(m,V,\frac{s}{4}) \neq \emptyset$ for almost
 every $m \in Y$.  Throwing away a set of measure zero, we may assume
 every point in $Y$ enjoys this property.  Next, Let $h(m) = \sup \{ d_M(m,y) | y \in Y \cap X(m,V,\frac{s}{4})\}$ when $m \in Y$.  Pick 
 $m_1 \in Y \cap X(m,V,\frac{s}{4})$ so that $d_M(m,m_1) \ge \frac{3
 h(m)}{4}$.  Further, we define
\[ C_m=Q_V^{-1}\left (Q_V\left (B_M\left (m,\frac{1}{c_0}\left (\frac{s
 h(m)}{4}\right )^{depth(M)}\right )\right )\right ) \]

We claim that $Y \cap C_m \subset X(m,2h(m),V,s) \cup
X(m_1,2h(m),V,s)$ for each $m \in Y$.  To prove this, suppose $y \in Y
\cap C_m$.  Since $Q_V$ is a
 H$\ddot{\text{o}}$lder map and our choice of $C_m$, we have that
 $d_M(Q_V(m),Q_V(y)) \le c_0 d_M(m,y)^\frac{1}{depth(M)} \le \frac{s
   h(m)}{4}$.  Now, if $h(m) < d_M(m,y)$ then the previous inequality
 and the fact the $y \in X(m,V,\frac{s}{4})$ imply that $d_M(m,y) \le
 h(m)$.  Thus, $d_M(m,y) \le h(m)$ and so $d_M(m_1,y) \le 2 h(m)$.
 If, further, $y \nin X(m_1,2 h(m), V, s)$ then, 
\begin{equation*}
\begin{split}
s d_M(m_1,y) &\le d_M(Q_V(m_1),Q_V(y)) \le  d_M(Q_V(m_1),Q_V(m)) +
d_M(Q_V(m),Q_V(y))\\
&< \frac{s}{4}d_M(m,m_1) + \frac{s h(m)}{4} \le \frac{s h(m)}{2}
\end{split}
\end{equation*}

Putting this together with the fact that $d_M(m,m_1) > \frac{s
  h(m)}{4}$, we have $d_M(m,y) > \frac{h(m)}{4} \ge
  \frac{d_M(Q_V(m),Q_V(y))}{s}$.  In other words, $y \in X(m,2
  h(m),V,s)$.  

From the hypothesis of the lemma, this gives 
\[ \Ha^k_M(Y \cap C_m) \le 2 \lambda (2 h(m)s)^{k depth(M)} \]

Using a standard Vitali covering argument, we can find a countable set
$I \subset Y$ s.t. $\{Q_V \left( B_M\left (y,\frac{1}{c_0} \left( \frac{s h(y)}{20}
\right )^{depth(M)}\right)\right)\}_{y \in I}$  are disjoint and $Q_V(Y)$
is contained in the union of the projected balls with the same centers
and five times the radius.  Therefore $Y \subset \cup_{y \in I} C_y$
and
\begin{equation*}
\begin{split}
\Ha^k_M(Y) &\le \sum_{y \in I} \Ha^k_M(Y \cap C_y) \le 2 \lambda 2^{k depth(M)}
\sum_{y \in I} (h(y)s)^{k depth(N)} \\
&= \frac{2^{k depth(M) +1}\lambda 20^{k depth(M)}}{\Ha^k_M(B_M(e^0,1) \cap
  exp(\mathfrak{n}))} \sum_{y \in I} \Ha^k_M\left ( exp\left (\mathfrak{n}\right) \cap
B_M\left( Q_V(y), \left (\frac{s h(y)}{20} \right )^{depth(M)}\right)\right)\\
& \le \frac{2^{k depth(M)+1}\lambda 20^{k depth(M)}}{\Ha^k_M\left(B_M(e^0,1) \cap
  exp\left(\mathfrak{n}\right)\right)} \Ha^k_M\left( exp(\mathfrak{n}) \cap
B_M\left(Q_V(m), \frac{\delta}{2}\right)\right) = C(k,M) \lambda \delta^k
\end{split}
\end{equation*}
$\qed$

 \begin{Cor} \label{mconecor}  Let $V, \mathfrak{n}$ and $s$ be as in the last
   lemma.  Suppose $Y \subset M$ is purely $N$-unrectifiable with
   $\Ha^k_M(Y) < \infty$, then there exists a constant $C$ depending
   only on $k$ and $\Ha^k_{N}(B_{N}(e^0,1))$ such that 
 \[ \Ha^k_M(Y \cap X(y,V,s),y) \ge C(sr)^{k depth(M)} \]
 for $H^k_M$ a.e. $y \in Y$.
 \end{Cor}

 \pf Again, this is the same as the proof of corollary 15.16 in
 \cite{Mattila} with the obvious modifications depending on the
 constant in lemma \ref{meascone}.  The only detail which must be checked is the lower
 bound on the density of sets of positive measure at almost every
 point.  This is covered in lemma \ref{densities}.

Now, we can prove the following theorem.

\begin{Thm} Suppose $Y$ is a $\Ha^k_M$-measureable subset of a Carnot
  group $M$ and $T$ is a subspace of $\mathfrak{m}$ so that $exp(T)$
 is isomorphic to $N$.  Further suppose that for $\Ha^k_M$ a.e. $y \in
 Y$ $T$ is a strong
  approximate tangent cone at $y$.
Then, $Y$ is $N$-rectifiable.
\end{Thm}

\pf
We will use lemma \ref{meascone} to show this implication via
contradiction.  Thus, we assume $Y$ is purely $N$-unrectifiable and
we will show that it can possess an approximate tangent cone
isomorphic to $N$ satisfying the hypotheses almost
nowhere.  In this proof, we will identify $M$ with its lie algebra,
denoted, for now, as $\R^n$.  Let
$k_0$ be the topological dimension of $N$ and let $G(n,k_0)$ be the
space of $k_0$-planes in $\mathfrak{m}$ which are isomorphic to $\mathfrak{n}$.  Let $m \in \mathbb{N}$ (we will
adjust the choice of $m$ later in the proof).  First, define the CC
operator norm on $V,W \in G(n,k_0)$ by $||V-W||_{cc}= sup_{\{u |
  d_M(e^0,e^u) <1\}} d_M(P_V(e^u),P_W(e^u)$.  Second, since the
identity map $(M,d_m) \ra \R^n$ is Lipschitz and the
inverse is H\"older, we see that $G(n,k_0)$ is compact with
respect to the CC operator norm as well as the standard operator norm.
We can use the 
compactness of $G(n,k_0)$ with respect to the CC operator norm to
cover $G(n,k_0)$ by finitely many balls 
of radius $\frac{1}{2m}$.  Call this family of balls $\EuScript{B}$.
Consider now the subset $\EuScript{G}=\{B_1,...,B_l\}$ of $\EuScript{B}$ of balls
containing $k_0$-planes that, under the induced grading from $N$, are
isomorphic to $N$.  Let $G_{gr}(n,k_0) \subset G(n,k_0)$ be the
subset of $k_0$-planes isomorphic to $N$ under the induced grading.
Now, for each $B_i$, pick $W_i \in B_i \cap G_{gr}(n,k_0)$.  Then, by 
construction, for $v \in B_i  \cap G_{gr}(n,k_0)$, we know that $||P_V
- P_{W_i}||_{cc} < \frac{1}{m}$.  In other words, if $\alpha$ is any vector, we have
\[d_M(P_V(e^\alpha),P_{W_i}(e^\alpha)) \le C(m)d_M(e^0,e^\alpha)\]
where $C(m)=const \cdot (\frac{1}{m})$.  Here, we abuse notation and
assume $V$ and $W_i$ 
are exponentiated at the same base point.

Fixing $W_i$, let $B=\left \{y \in Y \bigg | \exists V \in B_i
\text{ s.t. } V \in apTan^{N}(Y,y) \text{ and $(+)$ holds} \right \}$.  We want to show that $B$ has
$\Ha^k_M$-measure zero.  Suppose, on the contrary, that $B$ has
positive measure.   Since, by assumption, we know that every point in
$B$ has a strong approximate tangent cone in $B_i$, given $\lambda >0$,
$\exists r_0 >0$ so that
\[ C =  \left \{b \in B \bigg | \sup_{0 <r<r_0} \frac{ \Ha^k_M(B \cap B_M(b,r)
  \minus X(b,V,C(m)))}{r^{k depth(N)}} < \lambda \right \}\]
has positive $\Ha^k_M$-measure.  

Next, we claim that, for sufficiently large $m \in \mathbb{N}$, 
\[C \cap X(b,r,W_i^{\perp},C(m)) \subset C \cap (B_M(b,r) \minus X(b,V,C(m)))\]

Geometrically, this is almost clear, but we shall prove it anyway.
Suppose the claim does not hold;  for every $m$, there exists $y \in
C \cap (X(b,r,W_i^{\perp},C(m)) \cap X(b,V,C(m)))$ and $y \neq b$.  From the
definitions of the ``X'' sets, we have

\begin{gather}
d_M(P_{W_i}(y),P_{W_i}(b)) < C(m) d_M(y,b) \tag{$i$} \label{a} \\
d_M(P_{V^{\perp}}(y),P_{V^\perp}(b)) < C(m) d_M(y,b) \tag{$ii$} \label{b} 
\end{gather}

From above, we know $d_M(P_{W_i}(y),P_V(y)) \le C(m) d_M(b,y)$.  Using the
triangle inequality and (\ref{a}), we have
\begin{gather}
d_M(P_V(y),P_V(b)) < 3 C(m) d_M(y,b) \tag{$iii$} \label{c}
\end{gather}

Now, as $m \ra \infty$, $C(m) \ra 0$ so, considering (\ref{b}) and
(\ref{c}), we can conclude that there exists an $0 < s <1$
sufficiently small such that
$y \in X(b,V,s) \cap X(b,V^\perp,s)$.  By the second portion of the
definition of strong approximate tangent cone, no
such $y$ can exist.  
Picking $m$ so that the claim is true, we can now finish the proof of
this implication using the following computation:

\begin{equation*}
\begin{split}
\Ha^k_M(C \cap X(b,r,W_i^\perp,C(m))) &\le \Ha^k_M(C \cap (B_M(b,r)
\minus X(b,V,C(m))))\\
&\le \lambda r^{k depth(M)}
\end{split}
\end{equation*}

Picking $\lambda$ sufficiently small, we violate corollary
\ref{mconecor}.  Thus, $C$ has measure zero as does $B$ proving the
theorem. 
$\qed$
\newline \newline

\noindent
{\em Remark:} 
We reiterate that the class of $N$-rectifiable manifolds, for a
  Carnot group $N$, has good local properties.  In particular, the
  tangent cone at each point (in this setting, it is easy to see using
  the $N$-approximability at a point that
  the unique approximate tangent cone will coincide with the tangent
  cone in the sense of Gromov) is isometric to $N$ in contrast to the
  examples cited earlier.  This should allow for a stronger 
  analysis of the local geometry of $N$-rectifiable smooth submanifolds.

\section{Level Sets of $C^1_N$ Functions}
We turn now to a slightly different examination.  Considering, for
example, smooth submanifolds of Carnot groups, one is immediately led to a
situation which is not covered in the previous discussion.  Namely,
the previous discussion considered only the case when the tangent
objects were themselves Carnot groups.  Considering submanifolds in
the three dimensional Heisenberg group, one sees that the
smooth tangent spaces (one possible local approximator) are not
isomorphic to Carnot groups but to subgroups of Carnot groups (see the
explicit example below).

To shed some light on this situation, this section investigates the properties of level sets of $C^1_N$
functions on a Carnot group, $N$.  The main reason for focusing on
this case is that it ties in nicely with the beautiful work of
Franchi, Serapioni and Serra Cassano (\cite{FSSC2}) who investigate
rectifiability from a different point of view.  In the Heisenberg
group,$\mathbb{H}$, they show that a surface, $S$, of bounded perimeter can be
decomposed as
\[ S = Z \bigcup \cup_{i=1}^\infty K_i\]
where $Z$ is a set of perimeter zero and $K_i$ are compact subsets of 
noncharacteristic level sets of $C^1_\mathbb{H}$ functions.  Thus, it
is natural to ask whether level sets of $C^1_N$
functions possess rectifiability properties of the sort discussed in
this paper near noncharacteristic points.

In modeling the
local behavior of such sets, we will need to use subgroups of the
Carnot group onto which the projection mapping is not Lipschitz.  In fact,
as in the Euclidean case, the subgroup which locally models
$f^{-1}(t)$ at a point $x$ is $ker (df_x)$ which, by Pansu's
differentiability theory, exists almost everywhere and is a subgroup
for a.e. $x$.  

We begin with an example illustrating that the process of constructing
Lipschitz maps is quite delicate. \\

\begin{Example}
Consider the mapping $f:H^3 \ra \R$ given by
$e^{aX+bY+cZ} \ra ((|a|^2+|b|^2)^2+|c|^2)^\frac{1}{4}$.  This is a
Lipschitz map (this is a 
quasi-norm on $H^3$, although it is slightly different than the
quasi-norm defined in section 2, and hence  is biLipschitz equivalent to the
Carnot-Carath\'eodory metric) and the
inverse image of $t \in \R$ is the boundary of $B_{qn}(n_0,t)$.
We will now calculate the kernels of the differential mappings at each
point.  First, we calculate $df_x(e^{aX + bY + cZ})$.  Fix $x=e^{\alpha
  X + \beta Y +  \gamma Z}$.  Since, \[xh_re^{aX + bY + cZ}=
e^{(\alpha+ra)X + (\beta+rb)Y +(\gamma+r^2c + \frac{r}{2}(\alpha b - a
  \beta))Z}\] we have
\begin{gather*}
\lim_{r \ra 0^{+}} \frac {((|\alpha+ra|^2 +|\beta+rb|^2)^2+ |\gamma+r^2c
    + \frac{r}{2}(\alpha b - a \beta)|^2)^\frac{1}{4}- ((|\alpha|^2 + |\beta|^2)^2
    +|\gamma|^2)^\frac{1}{4}}{r} \\
= \frac{\alpha^3a+\gamma \alpha b + \alpha \beta^2a - \gamma \beta a
    + \alpha^2 \beta b + \beta^3 b}{(\alpha^4+2 \alpha^2 \beta^2 +
    \beta^4 + \gamma^2)^\frac{3}{2}}
\end{gather*}

So, for this fixed $x$, we see that when $\gamma
\alpha+\alpha^2\beta + \beta^3 \neq 0$, 
\[ ker (df_x) = \left \{ e^{aX+bY+cZ} \bigg | b= \frac{\gamma \beta -
    \alpha^3 -  \alpha \beta^2}{\gamma \alpha + \alpha^2 \beta +
    \beta^3} a \right \} \] 
Note that these subgroups are all mutually isomorphic.  One should note
that the case when $\alpha=\beta=0$ is particularly interesting, for
at those points, the kernel is $\{ e^{aX+bY+cZ}\}$.  These are points
where Pansu's differential is zero and hence are characteristic points
of the level set.  As mentioned above, we will be concerned only with
neighborhoods of noncharacteristic points.  
\end{Example}

For this section, we make the following assumptions:

\begin{itemize}
\item $N$ is a Carnot group of (topological) dimension
$n+m$, where $n$ is the dimension of the horizontal directions.
\item $k$ denotes the Hausdorff dimension of $N$.
\item $f:N \ra \R$ is a $C^1_N$ function.
\item For a fixed $t$, the characteristic set of $f^{-1}(t)$ is $\Sigma_{t}$.
\end{itemize}  

The most natural candidate for a
Lipschitz map between $U \subset N$ and $f^{-1}(t)$ would be, as in
the Euclidean case, the inverse of the projection of $f^{-1}(t)$ to
$ker (df_x)$.  Unfortunately, this map is not always Lipschitz.  We show this
in a special case.  Consider $x=e^X$ in $f^{-1}(1)$.  Then, by the
above calculation, $ker(df_x)=\{e^{bY+cZ}\}$.  Let $P_x: f^{-1}(1) \ra
ker(df_x)$ be the map $P_x(e^{aX+bY+cZ})=e^{bY+cZ}$ where
$((|a|^2+|b|^2)^2+|c|^2)^\frac{1}{4}=1$.  Taking $n_1=e^{aX+\sqrt{1-a^4}Z}$ and
  $n_2=e^{aX+\epsilon Y + \sqrt{1-(a^2+\epsilon^2)^2}Z}$ for $\epsilon >0$
  small and $a$ close to $1$, we have
\begin{equation*}
\begin{split}
d_{H^3}(n_1,n_2) &= d_{H^3}(e^0,e^{\epsilon Y + (\sqrt{ 1 - (a^2 +
  \epsilon^2)^2} -\sqrt{1-a^4}-\frac{1}{2}a\epsilon)Z}) \\
& \ge C_1 \epsilon^{\frac{1}{2}} 
\end{split}
\end{equation*}
But, computing the distance between the projected points, we have
\begin{equation*}
\begin{split}
d_{H^3}(P_x(n_1),P_x(n_2)) &= d_{H^3}(e^0, e^{\epsilon Y - (\sqrt{ 1 - (a^2 +
  \epsilon^2)^2} -\sqrt{1-a^4})
  Z})\\
&\le C_2 \epsilon 
\end{split}
\end{equation*}
Taking the quotient and applying the above estimates yields that
$Lip(P_x^{-1}) \ge C_3 \epsilon^{-\frac{1}{2}}$.  Letting $\epsilon$ go to
zero shows that the inverse map cannot be Lipschitz.

In the theory of $\R^n$-rectifiability, one has a nice decomposition
of Lipschitz maps between Euclidean spaces controlled by the area and
co-area formulae.  In particular, one can conclude that the inverse
image of a point under a Lipschitz map is $\R^k$-rectifiable.  We
prove, in this section, that level sets of $C^1_N$  maps from $N$ to
$\R$ are $T_x$-approximable at $x$, a nonchracteristic point, and have approximate tangent cones isomorphic to
$T_x$ at $x$ for $T_x$ some subgroup of $N$.  The proof of such a fact is significantly harder
than in the Euclidean case because many of the tools (even the ones
used in the previous sections) are not available in the cases where
the local model for the level sets are subgroups of Carnot groups with
incompatible gradings.  To circumvent the standard arguments, we
employ an approximation of the $C^1_N$ function by a smooth function
using the usual mollification procedure.  This gives us two tools -
the smooth approximation and the Pansu differential of the function.
Combining these two different perspectives allows us the desired
control.    

\begin{Lem} \label{smoothapprox}  If $f: N \ra \R$ is a $C^1_N$ map,
  there exist smooth functions $f_i:N \ra \R$ which converge uniformly
  to $f$ on compact sets.  In addition, if $X$ is a horizontal vector
  field on $N$, then $Xf_i \ra Xf$ uniformly on compact sets as
  $Xf$ is continuous.
\end{Lem}

\noindent
{\em Remarks on proof:}  This is a simple consequence of the
convergence properties of 
convolution operators on graded nilpotent Lie groups and the usual
mollification and approximation procedures.  As in the Euclidean case,
one has an approximation continuous functions by smooth functions
which converge uniformly on compact sets.  These facts are well known
and can be found in, for example, \cite{Hardy}.   The reader
should also see \cite{GarNh},
\cite{FSSC} and \cite{FSSC2} (in particular proposition 5.8 and theorem
6.4) for an explicit discussion of
mollifiers acting on functions on Carnot-Carath\'eodory spaces. 

% \begin{Cor}  If $T_{x_i}^i$ is the (smooth) tangent space to $f^{-1}_i(t)$
%   at a point $x_i$ then $T_{x_i}^i \ra T_x=ker(df_x)$ when $x_i \ra x$ and
%   $x$ is a point of
%   Pansu differentiability.
% \end{Cor}
 
% \pf This is a direct consequence of the second half of the previous
% lemma.  $\qed$

\begin{Lem} For almost every $t \in \R$, $f_i^{-1}(t)$ are all smooth
  submanifolds of $N$.
\end{Lem} 

\pf This is an application of Sard's theorem.

\begin{Lem} \label{diffpts} For $\Ha^k_N$ a.e. $t \in \R$, almost every point of $f_i^{-1}(t)$
  is a point of Pansu differentiability for $f$.  For almost every $t \in \R$,
  $f^{-1}(t)$ is $\Ha^{k-1}_N$ measurable and, for such a $t$, almost
  every $x \in f^{-1}(t)$ is a point of Pansu differentiability for $f$.
\end{Lem}
\pf  The first statement is a consequence of the coarea formula for smooth real
valued maps on Carnot groups.  See proposition \ref{coarea} and \cite{CCcalc} section
4.7.  The second statement follows directly from lemma 2.10.15 in \cite{Fed}.$\qed$

\begin{Thm}\label{Tapprox} Let $N$ be a Carnot group of Hausdorff
  dimension $k$ and let $f:N \ra \R$ be a
  $C^1_N$ map.  Then, for a.e. $t \in \R$ and $\Ha^{k-1}_N$ a.e. $x \in
  Int(f^{-1}(t)\minus \Sigma_t)$, there exists a subgroup
  $T_x \subset N$ and $r_x >0$ such that $B_N(x,r_x) \cap
  (f^{-1}(t)\minus \Sigma_t)$ is
  $T_x$-approximable at $x$ and $f^{-1}(t)$ has a unique approximate
  tangent cone isomorphic to $T_x$ at $x$.
\end{Thm}

The proof of this theorem follows from the next lemmas and lemma
\ref{nbd}.

\begin{Lem} \label{cond1}  Suppose $x \in f^{-1}(t)$ is a point of
  Pansu differentiability for $f$.  Let $T_x=x \cdot ker(df_x)$.  If $x' \in f^{-1}(t) \cap B_N(x,s)$ then
  $d_N(x',T_x)=o(s)$.
\end{Lem}

\pf Consider a one parameter family of points $xh_se^{w_1(s)+w_2(s)}
\in f^{-1}(t)$ where $w_1(s) \in ker(df_x)$, $w_2(s)$ is perpendicular
to the
kernel and $d_N(e^0,e^{w_1(s)+w_2(s)})=1$ for all $s$.  Now, by the
  choice of these points $|f(xh_se^{w_1(s)+w_2(s)}) - f(x)|=0$.
  Moreover, using the Pansu differentiability at the point $x$,
  $|f(xh_se^{w_1(s)+w_2(s)})-( f(x) +sdf_x(e^{w_2(s)}))|=o(s)$ and
  so, using the left invariance and homothety in the target,  we have:
\[ |df_x(e^{w_2(s)})|=o(1) \]
Moreover, we notice that $df_x|_{N/ker(df_x)}$ is biLipschitz since
$d_N|_{N/ker(df_x)}$ and $d_\R|_{Im(df_x)}$ are simply nondegenerate
  metrics on $\R$ and $df_x$ is a Lipschitz group homomorphism
  respecting dilations.  Using this,   we have that 
\[ d_N(e^{w_2(s)},e^0)=o(1)\]
So, using the Campbell-Baker-Hausdorff formula and the fact that
$w_2(s)$ lies in a Lie subalgebra,
$d_N(xh_se^{w_1(s)+w_2(s)},x\cdot ker(df_x))=d_N(h_{s}e^{w_2(s)},e^0)=o(s)$. $\qed$

Note that this lemma is vacuous at a characteristic point, where $ker(df_x)=N$.

Recalling the definition of Ahlfors regularity (definition
\ref{Ahlfors}) and the fact that smooth hypersurfaces of $N$ have
Hausdorff dimension $k-1$ (lemma \ref{hyper}), we show the Ahlfors
regularity of smooth level sets.  We remark that a
much more general version of this lemma is proved by R. Monti and
D. Morbidelli in \cite{MM}.  We include the proof here for
completeness.

\begin{Lem}\label{regular} Suppose $S$ is the level set of a smooth
  realvalued function on $N$ and let $\Sigma$ be its characteristic
  locus.  If $x \in S \minus \Sigma$ then there exists $r >0$
  so that the metric space $(\overline{B_N(x,r)\cap S}, d_N)$ is
  Ahlfors regular of dimension $k-1$ with constants depending only on
  $x$.  In other words, there exists a constant $A$ so that for $y
  \in\overline{B_N(x,r)\cap S}$ and $r>s >0$,
\[ A^{-1}s^{k-1} \le \Ha^{k-1}_N(B_N(y,s) \cap (\overline{B_N(x,r)\cap
  S})) \le A s^{k-1}\]
  
\end{Lem}

\pf To set up notation, 
 we assume $S$ is a smooth hypersurface given as a level set of a
smooth function $f:N\to\R$, i.e. $S=f^{-1}(0)$. Let
$\{X_1,\dots,X_{n}\}$ be an orthonormal basis for the
horizontal subbundle of $N$. Let $\{x_i\}$ be coordinates
associated to $\{X_i\}$ and let $\{z_1 , \dots , z_m\}$ be coordinates in the
nonhorizontal directions. Since  $x\in S \minus \Sigma$, $\n_0f$
does not vanish at $x$. Then, we can find a horizontal vector field,
$V$, such that $Vf(x)\neq0$. Without loss of generality, we may
take $V=X_1$. Using the (standard) implicit function theorem,
we can find a smooth function $g(x_2,\dots,x_{n},z_1, \dots, z_m)$ such
that $f(g(x_2,\dots,z_m),x_2,\dots,z_m)=0$, for
points in some neighborhood of $x$. In other words, $S$ is
parameterized via $g$:
\begin{equation*}
\begin{split}
\ov g:  &U\sbq\R^{n+m-1}\lra S\sbq N \cong \R^{n+m}\\
&(x_2,\dots,x_{n},z_1, \dots, z_m)\mapsto
(g(x_2,\dots,z_m), x_2,\dots, z_m).
\end{split}
\end{equation*}

With this notation in place, we now estimate the $\Ha_N^{k-1}$ measure
of sets of the form $Box(y,r)\cap S$.

First, pick a closed box around $x$, $C_x=\ov{Box(x,r) \cap S}$ inside $\ov{g(U)}$ 
so that $\n_0f$ and $f_{x_1}$ do not vanish on $C_x$. Then, $\exists$
positive constants $K_1$, $K_2$ s.t.
\begin{equation*}
\begin{split}
&K_1\le \left |\frac{\n_0f}{f_{x_1}} \right
|\le K_2 \qquad\text{for }\ y\in C_x.
\end{split}
\end{equation*}
Now, for $y\in C_x$, consider  $Box(y,s)\cap S\sbq C_x$.
By proposition \ref{surfmeas}
$$
\Ha^{k-1}_N(Box(y,s)\cap S \cap C_x) = \int_{Box(y,r)\cap S\cap C_x}
\Bigm|\f{\n_0f}{\n f}\Bigm| dA
$$
where $dA$ is the Riemannian area element. 

%Using the ball-box
%theorem, there exists a (universal!) constant $C>0$ s.t.
%$\cbox(y,C^{-1}r)\sbq B_N(y,r)\sbq\cbox(y,Cr)$. So,
%\begin{equation}
%%\begin{align}
%\int_{\cbox(y,C^{-1}r)\cap S}\Bigm|\f{\n_0f}{\n f}\Bigm| dA \le \int_{B_N(y,r)\cap S}\Bigm|\f{\n_0f}{\n f}\Bigm| dA  \le
%\int_{\cbox(y,Cr)\cap S}\Bigm|\f{\n_0f}{\n f}\Bigm| dA\tag{$*$}
%%\end{align}
%\end{equation}
To compute $dA$, we need to compute the determinant of the metric matrix.  First, note that $\ov g_{x_i} = (g_{x_i},0,
\dots, 1, \dots, 0)$.  Taking the matrix $G = (<\ov g_{x_i},\ov
g_{x_j}>)$ where $\ip$ is the Riemannian inner product on $N$.  Note
that, in the coordinates we are using for $N$, $\ip$ is the standard
Euclidean inner product.  Thus, 
\[ G= \begin{pmatrix}  g_{x_2}^2+1 & g_{x_2}g_{x_3}  & \dots & g_{x_2}g_{z_m} \\
  g_{x_3}g_{x_2} & g_{x_3}^2+1 &  \dots & g_{x_3}g_{z_m}\\
 &  \dots & \dots & \\
g_{z_m}g_{x_2} & g_{z_m}g_{x_3}&\dots &g_{z_m}^2+1 
\end{pmatrix}
\]

Computing $dA$ yields;
$$
(g_{x_2}^2+ \dots + g_{z_m}^2 + 1)^{\frac{1}{2}} dx_1\dots dx_{2n} dz_1
\dots dz_m.
$$
Using the implicit function theorem, we calculate that
$g_{x_j}=-\f{f_{x_j}}{f_{x_1}}\ \text{ for } n  \ge j>1$ and
$g_{z_i}=-\f{f_{z_i}}{f_{x_1}}$ for $1 \le i \le m$. So,
$dA=\f1{|f_{x_1}|}\cdot |\n f | dx_2\dots dz_m$. Thus,
$K_1dx_2\dots dz_m\le \left | \f{\n_0 f}{\n f} \right | dA\le
K_2dx_2\dots dz_m$. To conclude the computation, we need to find a
description of $Box(y, s) \cap S \cap C_x$.  Since we know bounds on $| \n
f|$, it follows that there exists a constant $K_3$ (in fact, we
may take $K_3$ to be one) so that if 
\[R_1=\left [-\frac{1}{K_3}s,\frac{1}{K_3}s \right]^{n-1} \cross \dots \cross
\left [-\left (\frac{1}{K_3}s\right)^{d(m)},\left
    (\frac{1}{K_3}s\right )^{d(m)} \right ]\]
then,
\begin{equation*}
 Box(y, s)\cap S \cap C_x \subset  \left \{ (g(x_2, \dots ,
    z_m),x_2, \dots, z_m) \bigg | (x_2, \dots,  z_m) \in R_1 \right \} 
\end{equation*}
Applying equation (*), we have 
\begin{equation*}
\Ha^{k-1}_N(Box(y,s) \cap S \cap C_x)  \le {K_2}\int_{R_1}dx_2\dots dz_m = 2 K_2
\frac{s^{k-1}}{K_3^{k-1}} 
\end{equation*}

To estimate $\Ha^{k-1}_N(Box(y,s)\cap S\cap C_x)$ from below, we
notice that, in the worst case, $y$ lies on the boundary of
$Box(x,r)$.  Thus, there is at least a quarter of $Box(y,s)$ inside of
$Box(x,r)$.  Using the same reasoning as above, there exists a
constant $K_4$ so that if 

\[R_2=\left [-\frac{1}{K_4}s,\frac{1}{K_4}s \right ]^{n-1}\cross
\dots  \cross
\left [-\left (\frac{1}{K_4}s \right )^{d(m)},\left
    (\frac{1}{K_4}s\right )^{d(m)} \right ] \]
then,
\begin{equation*}
 \left \{ (g(x_2, \dots ,
    z_m),x_2, \dots, z_m) \bigg | (x_2, \dots,  z_m) \in R_2 \subset  Box(y, s) \cap S \right \} 
\end{equation*}
and a quarter of this set lies inside $Box(y,s) \cap S \cap C_x$.

Again, using equation (*), we have 
\[ \frac{K_1}{2K_4^{k-1}}s^{k-1} \le \Ha^{k-1}_N(Box(y,s)\cap S\cap C_x)\]

To summarize, this says that the metric space $(C_x,d_{qn})$ is an
Ahlfors regular space of dimension $k-1$ where the constant is $\max
\left \{\frac{2K_2}{K_3^{k-1}}, \frac{2K_4^{k-1}}{K_1}\right \}$.

Now by the Ball-Box theorem, we have that 
\begin{multline*}
\frac{K_1}{2(CK_4)^{k-1}}s^{k-1} \le \Ha^{k-1}_N(Box_N(y, C^{-1}s)\cap S \cap
C_x) \le \Ha^{k-1}_N(B_N(y,s) \cap S\cap C_x) \\
\le  \Ha^{k-1}_N(Box_N(y, Cs)\cap S \cap
C_x) \le \frac{2K_2 C^{k-1}}{K_3C^{k-1}} s^{k-1}
\end{multline*}
In other words $(C_x,d_N)$ is an Ahlfors regular space of dimension
$k-1$.  Again using the Ball-Box theorem and potentially adjusting
$r$, we see that $(\overline{B_N(x,r)\cap S},d_N)$ is Ahlfors regular
of dimension $k-1$ with constant $A$ depending on
$\{C,K_1,K_2,K_3,K_4\}$.  Since, $C$ is a universal constant (for $N$)
and the $K_i$ depend only on the point $x$, we have the desired result.
$\qed$

\noindent
{\em Remark: } Suppose $x$ is a noncharacteristic point in $f^{-1}(t)$ with the property
that there exists a neighorhood of $x$ with no characteristic points.
Consider the case when we have smooth functions
$f_i$ approximating a $C^1_N$ function $f$.  Then, using lemma \ref{cond1},
for sufficiently large $i$ (i.e. when the approximation is close
enough), we can pick the constants $K_3$ and $K_4$ uniformly in $i$.
Also, since $\n_0 f_i$ converges to $\n_0 f$, we can pick the
constants $K_1$ and $K_2$ uniformly as well.  \newline

\begin{Lem} \label{measjump}  Fix a compact set $\Omega \subset N$.
  Let $C=(f^{-1}(t)\minus \Sigma_t) \cap \Omega$ and $C_i=f^{-1}_i(t) \cap \Omega$.
  Then, there exists a constant $\kappa > 0$ such that for $x \in N$ and $r >0$,
\[\kappa \Ha^{k-1}_N(C \cap B_N(x,r)) \ge \overline{\lim}_{i \ra \infty}
  \Ha^{k-1}_N(C_i \cap B_N(x,r))\]
\end{Lem}

\pf  This is not a hard fact to prove, but we cite it as a special
case of lemma 8.35 in \cite{brokendreams}. 
First we note that by the previous lemma, $C_i \cap
B_N(x,r)$ is Ahlfors regular of dimension $k-1$.  Second, by the
remark after the previous lemma, the constants can be picked
universally for large enough $i$. 
 Thus, the hypotheses of
lemma 8.35 in \cite{brokendreams} are satisfied.  $\qed$  

In the next lemma, we use the notation $Int(C)$ to denote the interior
of a set $C$.

\begin{Lem}\label{cond2} Fix $\alpha >0$.  Then,  for $x \in Int(f^{-1}(t)
  \minus \Sigma_t)$, there exists $s_0 >0$,
  $\theta >0$ such that for $0 < s< s_0$ and $x' \in T_x \cap
  B_N(x,s)$ then,
    $\Ha^{k-1}_N(f^{-1}(t)\cap B_N(x',\alpha s)) \ge \theta s^{k-1}$.  
\end{Lem}

\pf  To determine $s_0$, we use lemma \ref{cond1}.  The estimate in
this lemma and the triangle inequality tell us that if $x'\in T_x \cap
B_N(x,s)$ then $d_N(x',f^{-1}(t)) = o(s)$.  Picking $s_0$ small enough
(this will depend on $\alpha$), we can guarantee that $B_N(x',
\frac{\alpha s}{2}) \cap f^{-1}(t) \neq \emptyset$ for $0<s<s_0$ and $x' \in T_x \cap
B_N(x,s)$.  Let $B_s^i=B_N(x',\alpha s) \cap f^{-1}_i(t)$.  We know
from the construction of the $f_i$ that the $B_s^i$  converge
to $B_s=B_N(x',\alpha s)\cap f^{-1}(t)$.  Therefore, by lemma
\ref{measjump}, $ \overline{\lim}_{i \ra
  \infty} \Ha^{k-1}_N(B_s^i) \le  \kappa \Ha^{k-1}_N(B_s)$.  Also, in
lemma \ref{regular}, we noted that $f^{-1}_i(t) \cap
\Omega$ is Ahlfors regular.  By our choice of $s_0$, we know that, for
sufficiently all large $i$, there are points $x_i \in f^{-1}_i(t) \cap
\Omega$ such that $B_N(x_i,\frac{\alpha s}{4}) \cap f_i^{-1}(t)
\subset B_s^i$.  By Ahlfors regularity, we have $\Ha^{k-1}_N(B^i_s)
\ge \frac{K \alpha^{k-1} s^{-k-1}}{4^{k-1}}$.  Moreover, the remark
after lemma \ref{regular} allows us to pick a constant $K$ independent
of $i$, yielding the desired
estimate.  $\qed$
\newline

\noindent
{\em Proof of theorem \ref{Tapprox}:}  As in lemma \ref{cond1}, we let
$T_x=x \cdot ker(df_x)$.  Next, for this candidate local structure, we must verify conditions (1)
and (2) in definition \ref{Napp}.  Condition (1) is verified by lemma
\ref{cond2}.  Note that the constant $r_x$ in the statement of theorem
is given by the constant $s_0$ in lemma \ref{cond2}.  Moreover, a consequence of lemma \ref{cond1}, there
exists $r_0 >0$ such that $(f^{-1}(t) \cap B_N(x,r))
\minus T_x(\alpha r) = \emptyset$ for $0 < r < r_0$. The proof the
existence of a unique approximate tangent cone isometric to $N$
follows just as in one implication of theorem \ref{Equiv}.  $\qed$
\newline

\noindent
{\em Remark:}  The techniques and computations in this section are
very similar in spirit to those used to prove the implicit function
theorem (on the Heisenberg group) in \cite{FSSC2}. The author is
grateful to B. Franchi, R. Serapioni, and F. Serra Cassano for
providing an early version of \cite{FSSC2} as well as valuable
conversations concerning their techniques.

Next, we note that if the map $f$ is suitably generic at $x$, then a
neighborhood of $x$ is $T_x$-approximable.  Let $\mathfrak{n}$ be the
Lie algebra of $N$ and suppose $\{X_1,\dots,X_{n}\}$ is a basis for
$\mathcal{V}$.  We define $O \subset N$ by
\[O=\{a_1X_1+ \dots + a_n X_n | a_i \neq 0 \text{\; for all $i$}\}\]

\begin{Def}  Fix $t \in \R$.  A $C^1_N$ map $f:N \ra \R$ is {\bf
    generic} at $x \in f^{-1}(t)$ if $\nabla_0f(x) \in O$.
\end{Def}

\begin{Lem}\label{nbd}  Suppose $\Omega_t$ is a compact subset of
  $f^{-1}(t) \minus \Sigma_t$ and suppose $f$ is generic at $x \in
  \Omega_t$.  Then, there exists $r_x >0$ so that for each $y \in B_N(x,r_x) \cap
  \Omega_t$, $ker(df_y)$ and $ker(df_x)$ are isomorphic as subgroups
  of $N$.
\end{Lem}
   
%\begin{Lem} \label{nbd}  Let $\Omega$ be a compact set in $N$ and let
%  $\Omega_t=\Omega \cap f^{-1}(t) \minus \Sigma_t$.  Then, for almost every $t$ and for almost every $x \in
%  \Omega_t$, there exists $r_x >0$ such that for a.e. $y \in B_N(x,r_x) \cap
%  \Omega_t$, $ker(df_y)$ and $ker(df_x)$ are isomorphic as subgroups
%  of $N$.
%\end{Lem}

\pf  
%Pick $t$ so that
%$f_i^{-1}(t)$ are smooth, a.e. point of each $f_i^{-1}(t)$ is a point
%of Pansu differentiability for $f$ and a.e. point of
%$f^{-1}(t)$ is a point of Pansu differentiability for $f$.  
The exponential preimage of $ker(df_x)$ in the Lie algebra based
at $x$ is a vector subspace of $\mathfrak{n}$.  Denote this subspace
by $V_x$.  Recall that the Lie algebra $\mathfrak{n}$ is graded, i.e. 
\[\mathfrak{n}=n_1 \oplus n_2 \oplus \dots \oplus
n_l\]

We define the induced grading data on $V_x$ by the decomposition
\[ V_x = V_x\cap n_1 \oplus V_x\cap n_2 \oplus \dots \oplus
V_x\cap n_l\]
We denote  $V_x\cap n_i$ by $V_x^i$.  We say that two
subspaces, $V_x$ and $V_y$ have compatible induced grading data if there is a algebraic
isomorphism between $V_x$ and $V_y$ mapping $V_x^i$ to $V_y^i$ for all $i$.

Now, the subgroups $ker(df_x)$ and $ker(df_y)$ are
isomorphic and biLipschitz with respect to the quasi-norm of $N$
restricted to the respective subgroups if $V_x$ and $V_y$ have
compatible induced grading data.  In this case, we say that
$ker(df_x)$ is equivalent to $ker(df_y)$.  By the assumption of genericity and the continuity
of $df_x=\nabla_0f(x)$, we see that there exists an $r_x >0$ so that
$\nabla_0f(y) \in O$ for all $y \in B_N(x,r_x)\cap f^{-1}(t)$.  
Now, fixing a $y \in  B_N(x,r_x)\cap f^{-1}(t)$ we can find constants
so that 
\[V_x = span \{a_2X_1-X_2, \dots, a_n X_1 - X_n\}\]
\[V_y = span \{b_2X_1-X_2, \dots, b_n X_1 - X_n\}\]
for nonzero constants $\{a_i,b_i\}$.  We realize the compatibility
between $V_x$ and $V_y$ via the map which takes $a_iX_1-X_i$ to
$b_iX_1-X_i$.  Thus $ker(df_x)$ and $ker(df_y)$  are equivalent subgroups of $N$. $\qed$

\begin{Cor} Using the notation above, if $f$ is generic at $x \in
  f^{-1}(t)$ then there exists a neighborhood, $E$, of $x$ in $f^{-1}(t)$ and
  a subgroup $T_x$ of $N$ so that
  $E$ is $T_x$-approximable.  
\end{Cor}

\pf This follows directly from the previous lemma and the previous
theorem. $\qed$

\end{document}